%% file: Menger_and_consonant_sets_in_the_Sacks_model.tex
\documentclass[12pt]{amsart}
\usepackage{amsmath, amssymb, latexsym, amsthm, mathrsfs, color, mathtools, tikz,pgfplots,tikz-cd, graphicx, stackengine,url,accents, scrextend}
\usepackage[top=2.5cm, bottom=2.5cm, left=2.5cm, right=2.5cm]{geometry}

\input{shortcuts}

\usepackage[all]{xy}
\usepackage{float}
\title{Menger and consonant sets in the Sacks model}

		
\thanks{The research of the first and the third authors
was funded in whole by the Austrian Science Fund
(FWF) [10.55776/I5930 and 10.55776/PAT5730424].
The research of the second author was funded by the National Science Center, Poland Weave-UNISONO call in the Weave programme
Project: Set-theoretic aspects of topological selections 2021/03/Y/ST1/00122
}

\author[V.~Haberl]{Valentin Haberl}
\address{Institut f\"ur Diskrete Mathematik und Geometrie, Technische Universit\"at Wien, Wiedner Hauptstrasse 8-10/104, 1040 Wien, Austria.}
\email{valentin.haberl.math@gmail.com}
\urladdr{https://www.tuwien.at/mg/valentin-haberl/}
\author[P.~Szewczak]{Piotr Szewczak}
\address{Institute of Mathematics, Faculty of Mathematics and Natural Science,
College of Sciences, Cardinal Stefan Wyszy\'nski University in Warsaw, W\'oycickiego 1/3,
01–938 Warsaw, Poland}
\email{p.szewczak@wp.pl}
\urladdr{http://piotrszewczak.pl}
\author[L.~Zdomskyy]{Lyubomyr Zdomskyy}
\address{Institut f\"ur Diskrete Mathematik und Geometrie, Technische Universit\"at Wien, Wiedner Hauptstrasse 8-10/104, 1040 Wien, Austria.}
\email{lzdomsky@gmail.com}
\urladdr{https://dmg.tuwien.ac.at/zdomskyy/}
\begin{document}

\maketitle

\begin{abstract}
Using iterated Sacks forcing and topological games, we prove that the existence of a totally imperfect Menger set in the Cantor cube with cardinality continuum is independent from ZFC.
We also analyze the structure of Hurewicz and consonant subsets of the Cantor cube in the Sacks model.
\end{abstract}

\section{Introduction}

By \emph{space} we mean an infinite topological Tychonoff space.
A space $X$ is \emph{Menger} if for any sequence $\eseq{\cU}$ of open covers of $X$, there are finite families $\cF_0\sub\cU_0, \cF_1\sub\cU_1,\dotsc$ such that the family $\Un_{n\in\w}\cF_n$ covers $X$.
Every $\sigma$-compact space is Menger and every Menger space is Lindel\"of.
The Menger conjecture asserts that every subset of the real line with the Menger property is $\sigma$-compact.
By a result of Fremlin and Miller~\cite[Theorem~4]{FrMill}, this conjecture is false in ZFC. This opened a wide stream of investigations in the realm of special subsets of the real line.
The Menger property is a subject of research in the combinatorial covering properties theory but also appears in other branches of mathematics as local properties of function spaces~\cite{coc2}, forcing theory~\cite{ChoHru???, ChoRepZdo15} or additive Ramsey theory in algebra~\cite{Tsa18,color}.

The Menger property is closely related to infinite combinatorics.
Let $a,b\in\NN$.
We write $a\les b$ if the set $\sset{n}{a(n)>b(n)}$ is finite.
In such a case we say that the function $a$ is dominated by the function $b$.
A set $D\sub \NN$ is \emph{dominating} if any function in $\NN$ is dominated by some function from $D$.
Let $\fd$ be the minimal cardinality of a dominating subset of $\NN$.
The Menger property can be characterized in terms of continuous images, as follows: 
a set $X\sub\Cantor$ is Menger if and only if no continuous image of $X$ into $\NN$ is dominating.
This characterization was proved by Hurewicz~\cite[\S~5]{Hu27} and then much later but independently by Rec{\l}aw~\cite[Proposition~3]{reclaw}, so we call it the \emph{Hurewicz-Rec{\l}aw characterization of the Menger property}.
It follows that any subset of $2^\w$ with cardinality smaller than $\fd$ is Menger and there is a non-Menger set of cardinality $\fd$.

The above mentioned result of Fremlin and Miller is dichotomic, i.e., it splits ZFC into two cases using undecidable statements.
Bartoszy\'{n}ski and Tsaban provided in \cite{BaTs} a uniform ZFC counterexample to the Menger conjecture.
By \emph{set} with a given topological property we mean a space homeomorphic with a subspace of $\Cantor$.
A set is \emph{totally imperfect} if it does not contain a homeomorphic copy of $\Cantor$.

\bthm[{Bartoszy\'{n}ski, Tsaban~\cite[Theorem~16]{BaTs}}]\label{thm:bats}
There is a totally imperfect Menger set of cardinality $\fd$.
\ethm

In the first part of the paper we consider the following problem.

\bprb\label{prb:main}
Is there a totally imperfect Menger set of cardinality $\fc$?
\eprb
By Theorem~\ref{thm:bats}, it suffices to consider the case $\fd<\fc$.
In Section~2, we introduce a game characterization of the Menger property, which is one of the main tools needed in Section~3. We show that adding iteratively $\w_2$ Sacks reals to a ground model satisfying CH, we get a model where $\fd<\fc$ and the answer to Problem~\ref{prb:main} is negative.
In Section~\ref{Hech} we also prove that $\fd<\fc$ is consistent with the existence of  a totally imperfect Menger set of cardinality $\fc$.

In the second part of the paper we analyze the structure of \emph{Hurewicz} and \emph{consonant} subsets of  $\Cantor$.
Recall that a  
 space $X$ is \emph{Hurewicz} if for any sequence $\eseq{\cU}$ of open covers of $X$, there are finite families $\cF_0\sub\cU_0, \cF_1\sub\cU_1,\dotsc$ such that the family $\{\Un\cF_n: n\in\w\}$ is a \emph{$\gamma$-cover} of  $X$,
 i.e., the sets $\sset{n}{x\in \Un\cF_n}$ are cofinite for all  $x\in X$.
 Obviously, every $\sigma$-compact space is Hurewicz and every Hurewicz space is Menger.
Similarly to  the Menger property,
the Hurewicz property can be characterized in terms of continuous images, as follows.
A set $A\sub\NN$ is \emph{bounded} if there is a function $b\in\NN$ such that $a\les b$ for all $a\in A$.
A set $X\sub\Cantor$ is Hurewicz if and only if every continuous image of $X$ into $\NN$ is bounded.
This characterization was proved independently by Hurewicz~\cite[\S~5]{Hu27} and Rec{\l}aw~\cite[Proposition~1]{reclaw}, so we again call it the \emph{Hurewicz-Rec{\l}aw characterization of the Hurewicz property}.

\emph{Consonant} spaces were introduced by Dolecki, Greco and Lechicki~\cite{DolGreLec95} and characterized by Jordan~\cite[Corollary 11]{Jor20} in the following way.
Let $X\sub \Cantor$.
A cover of $X$ is a $\emph{k-cover}$ if any compact subset of $X$ is contained in some set from the cover.
A game $\GKO$ played on $X$ is a game for two players, \Alice{} and \Bob{}.
For a natural number $n$, in round $n$: \Alice{} picks an open $k$-cover $\cU_n$ of $X$ and \Bob{} picks a set $U_n\in\cU_n$.
\Bob{} wins the game if the family $\sset{U_n}{n\in\w}$ is a cover of $X$, and \Alice{} wins otherwise.
A set $Y\sub 2^\w$ is consonant if and only if \Alice{} has no winning strategy in the game $\GKO$ played on $2^\w\sm Y$.
We treat here this characterization as a definition of consonant sets.

Consonant spaces have close connections to combinatorial covering properties.
Let $Y\sub\Cantor$.
It follows from the game characterization of the Menger property given below that if the set $Y$ is consonant, then the set $\Cantor\sm Y$ is Menger.
A space $X$ is \emph{Rothberger} if for any sequence $\eseq{\cU}$ of open covers of $X$, there are sets $U_0\in\cU_0, U_1\in\cU_1,\dotsc$ such that the family $\sset{U_n}{n\in\w}$ is a cover of $X$.
Using a game characterization of the Rothberger property given by Pawlikowski~\cite{paw}, if the set $\Cantor\sm Y$ is totally imperfect, then $Y$ is consonant if and only if the set $\Cantor\sm Y$ is Rothberger. Indeed, if $2^\w\setminus Y$ is Rothberger, then Alice does not even have a winning
strategy in the game $\gone(\mathcal O,\mathcal O)$ on $2^\w\setminus Y$, hence
$ Y$ is consonant. Assuming now that $Y$ is consonant and $2^\w\setminus Y$ is totally imperfect, we shall show that  $2^\w\setminus Y$ is Rothberger.
Let $\mathcal U_{n,m}$ be an open cover of $2^\w\setminus Y$ for all $\langle n,m\rangle\in\w\times\w$. Set $\mathcal W_n=\{\bigcup_{m\in\w}U_{n,m}: U_{n,m}\in \mathcal U_{n,m} $ for all $m\in\w \}$. Since each compact subset of  $2^\w\setminus Y$ is countable,
each $\mathcal W_n$ is a $k$-cover of  $2^\w\setminus Y$.
Since Alice has no winning strategy in the game
$\gone(\mathcal K,\mathcal O)$ on  $2^\w\setminus Y$, we conclude that
 $2^\w\setminus Y$ is $\sone(\mathcal K,\mathcal O)$, and hence there are $W_n\in\mathcal W_n$,
$n\in\w$, such that $\bigcup_{n\in\w}W_n=2^\w\setminus Y$.
For each $n,m$ let $U_{n,m}\in\mathcal U_{n,m}$ be such that
$W_n=\bigcup_{m\in\w}U_{n,m}$. Then $2^\w\setminus Y=\bigcup_{n,m\in\w}U_{n,m}$,
which shows that $2^\w\setminus Y$ is Rothberger.\footnote{We have learned this argument from Paul Gartside.}


In Section~5 we introduce a topological game which allows us to analyze the structure of Hurewicz and consonant spaces in $\Cantor$.
We show that in the Sacks model, mentioned above, each consonant (Hurewicz) set $X\sub\Cantor$ and its complement $\Cantor\sm X$ are unions of $\w_1$ compact sets. 
This approach allows us also to capture a result of Miller~\cite[\textsection5]{mill} that in this model, every \emph{perfectly meager} subset of $\Cantor$, i.e., a set which is meager in any perfect subset of $\Cantor$, has size at most $\w_1$.

The main tools used in our investigations are game characterizations of the considered properties.
Let $X$ be a space.
The \emph{Menger game} played on $X$ is a game for two players, \Alice{} and \Bob{}.
For a natural number $n$, in round $n$: \Alice{} picks an open cover $\cU_n$ of $X$ and \Bob{} picks a finite family $\cF_n\sub\cU_n$.
\Bob{} wins the game if the family $\Un_{n\in\w}\cF_n$ covers $X$, and \Alice{} wins otherwise.
For more details about this game, we refer to the works of Scheepers~\cite[Theorem~13]{coc1} or Tsaban and the second named author~\cite{Mgame}.
Similarly to the Menger property, the Hurewicz property can also be characterized using a topological game.
The \emph{Hurewicz game} played on $X$ is a game for two players, \Alice{} and \Bob{}.
For a natural number $n$, in round $n$: \Alice{} picks an open cover $\cU_n$ of $X$ and \Bob{} picks a finite family $\cF_n\sub\cU_n$.
\Bob{} wins the game if the family $\sset{\Un\cF_n}{n\in\w}$ is a $\gamma$-cover of $X$, and \Alice{} wins otherwise.

\brem
In the definition of the Menger and Hurewicz games we could assume, in addition, that none of the
$\mathcal U_n$'s contains a finite subcover, and get an equivalent definition. Indeed,
if $X$ is compact, then Alice has no legal moves, and we standardly adopt the convention
that the player having no moves loses. Thus, Bob has a ``trivial'' winning strategy.
On the other hand, if $X$ is not compact, then it has a cover $\mathcal U_*$ without finite subcovers,
and there is no loss of generality in assuming that Alice always plays refinements of $\mathcal U_*$.
Similar comments apply  to the definitions of the Menger and Hurewicz properties.
\erem

\bthm[Hurewicz] \label{thm:hur}
A set $X\sub\Cantor$ is Menger if and only if \Alice{} has no winning strategy in the Menger game played on $X$.
\ethm

\bthm[{Scheepers~\cite[Theorem~27]{coc1}}] \label{thm:hur_hur}
A set $X\sub\Cantor$ is Hurewicz if and only if \Alice{} has no winning strategy in the Hurewicz game played on $X$.
\ethm

\section{Menger game and perfect sets} \label{sec2}

Now we shall address some specific instances of perfect spaces and families of their clopen subsets.
Suppose that $\seq{F_n}{n\in\w}$ is a  non-decreasing  sequence 
of finite subsets of some well-ordered set $\la S,<\ra$  such that $S=\bigcup_{n\in\w}F_n$, and $\seq{k_n}{n\in\w}$
is a strictly increasing sequence of natural numbers.
For each $n\in\w$ let $\Sigma_n\sub(2^{k_n})^{F_n}$.
Fix natural numbers $n,m$ with $n<m$.
For $\nu\in \Sigma_n$ and
$\sigma\in\Sigma_m$, we write $\nu\prec\sigma$ if $\nu$ is \emph{extended} by $\sigma$,
i.e., $\sigma(\beta)\restrict k_n=\nu(\beta)$ for all $\beta\in F_n$.
Let $C\sub\Sigma_n$.
A map $e\colon C\to\Sigma_m$ such that $\nu\prec e(\nu)$ for all $\nu\in C$, is \emph{coherent}, 
if  for any $\nu,\nu'\in C$, letting  $\beta\in F_{n}$ be
 the minimal element of $F_{n}$ with $\nu(\beta)\neq\nu'(\beta)$,
we have $e(\nu)\restrict (F_{m}\cap\beta)=e(\nu')\restrict (F_{m}\cap\beta)$.
In what follows we shall assume that
\begin{itemize}
\item[$(e_f)$] For every $C\subset\Sigma_n$, coherent $e_0:C\to \Sigma_m$, 
and $\nu\in\Sigma_n\setminus C$, 
there exist two different coherent maps $e,e'\colon (C\cup\{\nu\}) \to\Sigma_m$ such that
$e\restrict C=e'\restrict C=e_0$.
\end{itemize}
Note that $(e_f)$ applied to $C=\emptyset$ yields that 
for every $\nu\in\Sigma_n$ there exists at least two $\sigma\in\Sigma_m$
such that $\nu\prec\sigma.$
Moreover, using~$(e_f)$ iteratively for all natural numbers $n,m$ with $n<m$, set $C\sub\Sigma_n$, elements $\nu\in C$ and $\sigma\in \Sigma_m$ with $\nu\prec\sigma$, there is a coherent map $e\colon C\to \Sigma_m$ with $e(\nu)=\sigma$.

The objects defined above give rise 
to the perfect subset $K\sub (2^\w)^S$ consisting of those $x$ such that 
for every $n$ there exists $\nu\in\Sigma_n$ such that 
$x\in[\nu]$, where 
\[
[\nu]:=\sset{x\in (\Cantor)^S}{x(\beta)\restrict k_n=\nu(\beta)\text{ for all }\beta\in F_n}.
\]

For every $\beta\leq S$, let $\mathrm{pr}_\beta\colon K\to (2^\w)^\beta$
be the projection map.
For a set $C\sub \Sigma_n$, a map $E\colon C\to K$ is a
\emph{coherent selection} if for every $\nu\in C$ we have $E(\nu)\in[\nu]$ and for every $\nu,\nu'\in C$ and  $\beta\in F_{n}$, which is the minimal element in $F_{n}$ with $\nu(\beta)\neq\nu'(\beta)$, we have $\mathrm{pr}_\beta(E(\nu))=\mathrm{pr}_\beta(E(\nu'))$.

\blem \label{coh_sel}
In the notation above, if a set $X\sub K$ is totally imperfect, then for every $n$ there exists a coherent selection
$E\colon\Sigma_n\to K$ such that $E[\Sigma_n]\sub K\sm X$.
\elem

\bpf
Fix $n$ and enumerate $\Sigma_n$ injectively as $\{\nu_0,\ldots,\nu_N\}$.
Since $[\nu_0]\cap K$ is perfect by $(e_f)$, we can pick $E_0(\nu_0)\in ([\nu_0]\cap K)\setminus X$.
Fix a number $k<N$ and put $C:=\{\nu_i:i\leq k\}$.
Assume that we have already defined a coherent map
$E_k:C\to K\setminus X$.
Then for every $m>n$ and $\nu_i\in C$, there is $e^m_0(\nu_i)\in\Sigma_m$ such that $E_k(\nu_i)\in[e^m_0(\nu_i)]$.
In that way we define a map $e^m_0\colon C\to \Sigma_m$.
It is clear that this map is coherent.
By~$(e_f)$ there are $\mu_{\la 0\ra}\neq\mu_{\la 1\ra}\in\Sigma_{n+1}$
such that both 
$e^{n+1}_{\la 0\ra}=e^{n+1}_0\cup\{\la \nu_{k+1}, \mu_{\la 0\ra}\ra\}$
and $e^{n+1}_{\la 1\ra}=e^{n+1}_0\cup\{\la \nu_{k+1}, \mu_{\la 1\ra}\ra\}$
are coherent as maps from $C\cup\{\nu_{k+1}\}$ to $\Sigma_{n+1}$.

Suppose that for some $m>n$ we have defined 
mutually different $\{\mu_s:s\in 2^{m-n}\}\sub\Sigma_m\sm e_0^m[C]$
such that 
$e^{m}_{ s}=e^{m}_0\cup\{\la \nu_{k+1}, \mu_{ s}\ra\}$
is
coherent as a map from $C\cup\{\nu_{k+1}\}$ to $\Sigma_{m}$ for all 
$s\in 2^{m-n}$.
By~$(e_f)$ for every $s\in 2^{m-n}$ there are $\mu_{ s\hat{\ \ } 0}\neq\mu_{ s\hat{\ \ }1}\in\Sigma_{m+1}$
such that 
$$ \{\la e^m_0(\nu_i), e^{m+1}_0(\nu_i)\ra:i\leq k\}\cup \{\la \mu_s, \mu_{ s\hat{\ \ }j} \ra\} $$
is coherent as a map from $\Sigma_m$ to $\Sigma_{m+1}$ for all $j\in 2$.
It follows that
$$ e^{m+1}_{s\hat{\ \ }j}:=e^{m+1}_0\cup\{\la \nu_{k+1}, \mu_{ s\hat{\ \ }j}\ra\}$$
is
coherent as a map from $C\cup\{\nu_{k+1}\}\to\Sigma_{m+1}$ for all 
$s\in 2^{m-n}$ and $j\in s$, which completes our construction
of the maps $e^m_s$ and elements $\mu_s$ as above for \emph{all}
$m>n$ and $s\in 2^{m-n}$.

For every $t\in 2^\w$ let $z_t\in K$ be the unique element in $\bigcap_{m>n}[\mu_
{t\restrict (m-n)}]$ and note that $z_t\neq z_{t'}$ for any $t\neq t'$ in $2^\w$.
Since the set $\sset{z_t}{t\in 2^\w}$ is perfect, there exists $t$ with $z_t\not\in X$.
It suffices to 
observe that 
$E_{k+1}:=E_k\cup\{\la\nu_{k+1},z_t\ra\}$ is a coherent map whose range is disjoint from 
$X$, which allows us to complete our proof by induction on $k<N$.
\epf

\blem \label{clm:easy}
In the notation used above, let 
$X\sub K$ be a totally imperfect Menger set.
Then there exists a sequence 
$\seq{\la i_n,j_n,C_n\ra}{n\in\w}$  such that
\begin{enumerate}
\item $\seq{i_n}{n\in\w}$ is a strictly increasing sequence of natural numbers;
\item $C_n\sub\Sigma_{i_n}$; 
\item for every $n\in\w$, we have $j_n\in [i_n,i_{n+1})$,
and for every $\nu\in C_n$ there exists $e_n(\nu)\in \Sigma_{j_n}$ extending
$\nu$, such that 
$$C_{n+1}=\bigcup_{\nu\in C_n}\{\sigma\in\Sigma_{i_{n+1}}\: :\: \sigma \succ e_n(\nu)\};$$
\item the maps $e_n\colon C_n\to \Sigma_{j_n}$ are coherent; and
\item $\bigcap_{n\in\w}\bigcup_{\nu\in C_n}[\nu]\cap X=\emptyset$
\end{enumerate}
\elem
\bpf
We shall describe a strategy $\S$ for \Alice{} in the  Menger game played on $X$  such that
each play lost by \Alice{} gives rise to the objects whose existence we need to establish. For every $n\in\w$, let $E_n\colon\Sigma_n\to K$ be a coherent selection from Lemma~\ref{coh_sel}. 
Put 
\[i_0:=0,\quad C_0:=\Sigma_{i_0},\quad Z_0:=E_{i_0}[\Sigma_{i_0}]\sub K\setminus X.
\]
For every $j\geq i_0$ and $\nu\in\Sigma_{i_0}$, let $
\sigma_{0,j}(\nu)$ be the unique element of $\Sigma_j$
such that $E_{i_0}(\nu)\in [\sigma_{0,j}(\nu)]$.
Then,
\[
\smallmedset{\Un_{\nu\in\Sigma_{i_0}}[\sigma_{0,j}(\nu)]}{ j\geq i_0}
\]
is a decreasing family of clopen sets in $(2^\w)^S$, whose intersection is equal to the set $Z_0$.
Since $Z_0\sub K\sm X$, the family $\cU_0$ of all sets
\[
U^0_j:=K\setminus \bigcup_{\nu\in\Sigma_{i_0}}[\sigma_{0,j}(\nu)],
\]
where $j\geq i_0$, is an increasing open cover of $X$.
Now, $\S$ instructs \Alice{} to start the play with $\cU_0$.

Suppose that \Bob{} chooses 
$U^0_{j_0}$ for some $j_0\geq i_0$.
For each $\nu\in\Sigma_{i_0}$, let $e_0(\nu):=\sigma_{0,j_0}(\nu)$, an element of $\Sigma_{j_0}$.
Since $E_{i_0}$ is a coherent selection, the map $e_0\colon C_0\to \Sigma_{j_0}$ is coherent.

Then we put 
\[
i_1:=j_0+1,\quad  C_{1}:=\Un_{\nu\in C_0}\sset{\sigma\in\Sigma_{i_{1}}}{\sigma \succ e_0(\nu)}.
\]

Suppose that a natural number $i_n$ and a set $C_n\sub \Sigma_{i_n}$ have already been defined for 
some $n>0$.
Let
$Z_n:=E_{i_n}[C_{n}]\sub K\setminus X$.
For every $j\geq i_n$ and $\nu\in C_n$, let $
\sigma_{n,j}(\nu)$ be the unique element of $\Sigma_j$
such that $E_{i_n}(\nu)\in [\sigma_{n,j}(\nu)]$.
Then,
\[
\smallmedset{\Un_{\nu\in C_{n}}[\sigma_{n,j}(\nu)]}{ j\geq i_n}
\]
is a decreasing family of clopen sets in $(2^\w)^S$, whose intersection is equal to the set $Z_n$.
Since $Z_n\sub K\sm X$, the family $\cU_0$ of all sets
\[
U^n_j:=K\setminus \bigcup_{\nu\in C_{n}}[\sigma_{n,j}(\nu)],
\]
where $j\geq i_n$, is an increasing open cover of $X$.
Now, $\S$ instructs \Alice{} to play the family $\cU_n$.

Suppose that \Bob{} chooses 
$U^n_{j_n}$ for some $j_n\geq i_n$. 
Since $E_{i_n}$ is a coherent selection, the map $e_n\colon C_n\to \Sigma_{j_n}$ is coherent.
Then we put 
\[
i_{n+1}:=j_n+1,\quad  C_{n+1}:=\Un_{\nu\in C_n}\sset{\sigma\in\Sigma_{i_{n+1}}}{\sigma \succ e_n(\nu)}.
\]

This completes our definition of the strategy 
$\S$ for \Alice{} such that each infinite play in the Menger game on
$X$ in which \Alice{} uses $\S$ gives rise to a
 sequence
\[
\seq{\la i_n,j_n,C_n,e_n,\U_n\ra}{n\in\w}
\]
as described above. In particular, conditions
$(1)$--$(4)$ are satisfied by the construction.
By Theorem~\ref{thm:hur}, there is a play, where \Alice{} uses the strategy $\S$ and the play is won by \Bob{}.
Then
$X\sub\bigcup_{n\in\w}U^n_{j_n}$, i.e., 
\[
\emptyset=X\cap \bigcap_{n\in\w}\bigcup_{\nu\in C_n}[\sigma_{n,j_n}(\nu)]
\]
For each $n$, we have 
\[
\Un_{\nu\in C_n}[\sigma_{n,j_n}(\nu)]=\Un_{\nu\in C_n}[e_n(\nu)]=\Un_{\nu\in C_n}\bigcup\sset{[\sigma]}{\sigma\in\Sigma_{i_{n+1}},\sigma \succ e_n(\nu)}=\bigcup_{\sigma\in C_{n+1}}[\sigma],
\]
and thus

\[
\emptyset=X\cap \bigcap_{n\in\w}\bigcup_{\nu\in C_n}[\sigma_{n,j_n}(\nu)]=X\cap \bigcap_{n\in\w}\bigcup_{\nu\in C_{n+1}}[\nu].
\]
It follows that condition $(5)$ is also satisfied.
\epf


\section{Combinatorics of conditions in the iterated Sacks forcing} \label{sec33}

Here we deal with countable support iterations of the forcing notion introduced by Sacks~\cite{sacks}. We do not prove any essentially new results 
about these iterations in this section, but rather ``tailor''
 several results established in \cite{mill,BauLav} and perhaps somewhere else
 for
the purposes we have in Sections~\ref{menger_size} and \ref{menger_size_rev}.
We try to  follow notations used in \cite{BauLav}.

Let $2^{<\w}:=\Un_{n\in\w}2^n$.
For elements $s,t\in 2^{<\w}$, we write $s\sub t$ if the sequence $s$ is an initial segment of the sequence $t$, i.e., $s(i)=t(i)$ for all $i\in\dom(s)$.
A \emph{Sacks tree} is a set $p\sub 2^{<\w}$ such that for every $s\in p$ and a natural number $n$, we have $s \restriction n\in p$ and there are elements $t,u\in p$ with $s\sub t$, $s\sub u$, $t\not\sub u$ and $t\not\sub u$.
For Sacks trees $p$ and $q$, a condition $q$ is stronger than $p$ which we write $q\geq p$ if $q\sub p$. The \emph{Sacks poset} $\bbS$ is the set of all Sacks trees ordered by $\geq$.
For $p,q\in\bbS$ and natural numbers $m>n$, we write $(q,m)\geq(p,n)$ if $q\sub p$ and for every $s\in p\cap 2^n$, there are different elements $t,u\in q\cap 2^m$ such that $s\sub t$ and $s \sub u$.
For $p\in\bbS$ and $s\in p$, let $p_s:=\sset{t\in p}{t\sub s\text{ or }s\sub t}$.

Let $\alpha$ be an ordinal number and $\mathbb S_\alpha$ be an iterated forcing of length $\alpha$ with countable support, where each iterand is a Sacks poset.
For $p,q\in\mathbb S_\alpha$, let $q\geq p$ if $\supp(q)\supseteq \supp(p)$ and for every $\beta\in\supp(p)$, we have $q\restrict \beta\forces_\beta q(\beta)\geq p(\beta)$.

Let $p\in\mathbb S_\alpha$, $F\sub\alpha$ be  a finite set, $n$  a natural number and $\sigma\colon F\to 2^n$  a map.
If $F=\emptyset$, then the map $\sigma$ is \emph{consistent} with $p$ and $p|\sigma:=p$.
Assume that $\beta$ is the greatest element in $F$, the map $\sigma\restrict\beta$ is consistent with $p$ and the condition $p|(\sigma\restrict\beta)$ has already been defined.
If 
$$\big(p|(\sigma\restrict \beta)\big)\restrict \beta\forces_\beta \sigma(\beta)\in p(\beta),$$ then the map $\sigma$ is \emph{consistent} with $p$ and we define
$p|\sigma(\gamma)$ to be the following $\mathbb S_\gamma$-name $\tau$:
\begin{itemize}
\item $p|(\sigma\restrict \beta)(\gamma),$  if $\gamma<\beta$;
\item If $\gamma=\beta$, then  $\big(p|(\sigma\restrict \beta)\big)\restrict \beta\forces_\beta \tau=p(\beta)_{\sigma(\beta)}$ 
and $r\forces_\beta \tau=p(\beta)$ for $r\in\mathbb S_\beta$ incompatible with $\big(p|(\sigma\restrict \beta)\big)\restrict \beta$; and
\item $p(\gamma)$,  otherwise.
\end{itemize}

The following fact can be established by induction on 
$|F|$ in a rather straightforward way.
\bobs \label{obs:restr1}
In the notation above, if $\sigma$ is consistent with $p$,
then $\sigma\restrict\beta$ is consistent with both
$p\restrict\beta$ and $p$,  
and $\big(p|(\sigma\restrict \beta)\big)\restrict \beta=(p\restrict \beta)|(\sigma\restrict \beta)$.
\eobs

A condition $p\in\mathbb S_\alpha$ is \emph{$(F,n)$-determined}, where 
$F\sub \alpha$ and $n\in\w$, if every map $\sigma\colon F\to 2^n$ is either consistent with $p$, or there is $\beta\in F$ such that $\sigma\restrict \beta$ is consistent with $p$ and $(p\restrict \beta)|(\sigma\restrict \beta)\forces_\beta \sigma(\beta)\notin p(\beta)$.
For $q\in\mathbb S_\alpha$ and a natural number $m>n$, we write $(q,m)\geq_F(p,n)$ if $q\geq p$ and for every $\beta\in F$, we have $q\restrict \beta\forces_\beta (q(\beta),m)\geq(p(\beta),n)$.

Let $p$ be an $(F,n)$-determined condition.
We write $(q,n)\geq_F (p,n)$ if $q\geq p$ and every map $\sigma\colon F\to 2^n$ consistent with $p$, is also consistent with $q$. Next, we collect rather straightforward  facts about the notions introduced above, these are used in nearly all works investigating iterations of the Sacks forcing
or similar posets consisting of trees.

\bobs \label{obs:det2}
 Let $p\in\mathbb S_\alpha$ be an $(F,n)$-determined condition, $\Sigma$  the set of all maps $\nu\colon F\to 2^n$ consistent with $p$ and $\beta<\alpha$. 
Then the following assertions hold.
\begin{enumerate}
\item  If 
 $(q,n)\geq_F (p,n)$, then $q$ is also $(F,n)$-determined;
\item  If  $\sigma\in \Sigma$, then 
$p|(\sigma\restrict\beta)$ is $(F\setminus\beta, n)$-determined,
and $\nu\in (2^n)^{F\setminus\beta}$ is consistent with $p|(\sigma\restrict\beta)$ iff $(\sigma\restrict\beta)\cup\nu\in\Sigma$; 
\item The set $\sset{p|\sigma}{\sigma\in\Sigma}$ is a maximal antichain above $p$;
\item  $p$ is $(F\cap\beta,n)$-determined and $\{\sigma\restrict(F\cap\beta):\sigma\in\Sigma\}$
is the family of all functions from $F\cap\beta$ to $2^n$ consistent with $p$;
\item If $\sigma\in\Sigma$ and $r\geq p|\sigma$, then there exists 
$q\in\mathbb S_\alpha$ with $(q,n)\geq_F(p,n)$ and $q|\sigma=r;$
\item If $D\sub \mathbb S_\alpha$ is open and dense, then
there exists $q\in\mathbb S_\alpha$ with $(q,n)\geq_F(p,n)$ and $q|\sigma\in D$
for all $\sigma\in\Sigma$; and
\item If $\tau$ is an $\mathbb S_\alpha$-name for a real, 
then for each natural number $l$ there is a condition $q\in\mathbb S_\alpha$ and a family 
$\{y_\sigma:\sigma\in\Sigma\}\sub  2^l$
such that $(q,n)\geq_F (p,n)$ and $q|\sigma \forces \tau\restrict l= y_\sigma$
for all $\sigma\in\Sigma$.
\end{enumerate}
\eobs

The last three items of Observation~\ref{obs:det2} imply the following easy
 fact.

\blem[{Miller~\cite[Lemma~2]{mill}}]\label{lem:millstep1}
Let $p\in\mathbb S_\alpha$ be an $(F,n)$-determined condition and $\tau$  an $\mathbb S_\alpha$-name for a real such that $p\forces \tau\in \Cantor\sm V$.
Then for each finite set $Y\sub \Cantor\cap V$, there is a finite set $X\sub\Cantor\cap V$ disjoint from $Y$ such that $|X|\leq 2^{n\cdot|F|}$, and for each natural number $l$, there is a condition $q\in\mathbb S_\alpha$  such that $(q,n)\geq_F (p,n)$ and 
\[
q \forces \Exists{x\in X}\End{\tau\restrict l= x\restrict l}.
\] 
\elem

 Let  $\dot{G}_{\alpha}$ be an $\mathbb S_\alpha$-name for 
 a $\mathbb S_\alpha$-generic filter, and for
 $\beta<\alpha$ let $\mathbb S_{\beta,\alpha}$ be an
 $\mathbb S_\beta$-name for the iteration from 
 (including) $\beta$ to $\alpha$, so that 
 $\mathbb S_\alpha$ is forcing equivalent to $\mathbb S_\beta*\mathbb S_{\beta,\alpha}$. 
Whenever we work in the forcing extension $V[G_\beta]$ for some 
$\mathbb S_\beta$-generic filter $G_\beta$, we denote by
$\dot{G}_{\beta,\alpha}$ a $\mathbb S_{\beta,\alpha}^{G_\beta}$-name for a $\mathbb S_{\beta,\alpha}^{G_\beta}$-generic filter over $V[G_\beta]$. We shall need the following easy observation, we use the notation from the above.

\bobs \label{obs:cons3}
Suppose that $p$ is $(F,n)$-determined, $\beta\leq\alpha$, and $p\restrict\beta\in G_\beta$.
Then in $V[G_\beta]$, $p\restrict[\beta,\alpha)^{G_\beta}\in\mathbb S_{\beta,\alpha}^{G_\beta}$ is $(F\setminus \beta,n)$-determined. 

Moreover, if $\sigma\in (2^n)^F$ is consistent with $p$
and $(p|(\sigma\restrict\beta))\restrict \beta\in G_\beta$, 
then $\sigma\restrict(F\setminus\beta)$ is consistent with 
$p\restrict[\beta,\alpha)^{G_\beta}$ in $V[G_\beta]$; and if 
$(p|(\sigma\restrict\beta))\restrict\beta\in G_\beta$ and
$\nu\in (2^n)^{F\setminus\beta}$ is consistent with $p\restrict[\beta,\alpha)^{G_\beta}$ in $V[G_\beta]$, then $(\sigma\restrict\beta)\cup\nu\in (2^n)^F$
is consistent with $p$.
\eobs


\blem[{Miller~\cite[Lemma~5]{mill}}]\label{lem:stepmill}
Let  $p\in\mathbb S_\alpha$ be an $(F,n)$-determined condition and $\tau$ be an $\mathbb S_\alpha$-name for a real such that
\[
p\forces \tau\in \bigl(\Cantor\cap V[\dot{G}_{\alpha}]\bigr)\sm \Un_{\beta<\alpha}\bigl(\Cantor\cap V[\dot{G}_\beta]\bigr).
\]
Then for any $k\in\w$ there exist a condition $q\in\mathbb S_\alpha$, a natural number $l>k$, and elements $y_{\sigma}\in 2^{l}$, for all maps $\sigma\colon F\to 2^n$ consistent with $p$, with the following properties:
\be
\item $(q, n) \geq_F (p, n)$,
\item $q|\sigma\forces \tau\restrict l= y_{\sigma}$,
\item the maps $y_{\sigma}$ are pairwise different.
\ee 
\elem

\bpf
Let $\xi=\min(F)$ and note that the fact that $p$ is $(F,n)$-determined yields   $N\leq 2^n$
and $\{s_i:i<N\}\sub  2^n$ such that
$p\restrict\xi$ forces $p(\xi)\cap 2^n=\{s_i:i<N\}$.
For every $i<N$ let $\mu_i$ be the map  $\{\la \xi,s_i\ra\}$
and note that $\mu_i$ is consistent with $p$.

By induction on $i<N$,  using Lemma~\ref{lem:millstep1} and  Observation~\ref{obs:det2}(2), we can find  mutually disjoint finite sets $X_i\sub\Cantor\cap V$,  $i<N$, such that  $|X_i|\leq 2^{n\cdot (|F|-1)}$, and for each natural number $l$ there is a condition $u^l_i\in\mathbb S_\alpha$ with 
$u^l_j\restrict\xi\geq u^l_i\restrict\xi$ for all $i<j\leq N$, $u^l_0\restrict\xi\geq p\restrict\xi$, $(u^l_i,n)\geq_{F\setminus\{\xi\}} (p|\mu_i,n)$ and 
\[
u^l_{i}\forces \Exists{x\in X_i}\End{\tau\restrict l= x\restrict l}.
\]
Pick a natural number $l_{*}>k$ such that $x\restrict l_*\neq x'\restrict l_*$ for any distinct $x,x'$ in $\bigcup_{i<N} X_i$. As a result, the elements of the family $\sset{X_i\restrict l_*}{i<N}$ are mutually disjoint.

Now we proceed by induction on the cardinality of $F$.
If  $F=\{\xi\}$, then
 $|X_i|=1$ for all $i<N$ (because $ 2^{n\cdot (|F|-1)}=1$), i.e., 
$X_i=\{x_i\}$ for some $x_i\in 2^\w$.
Put $y_{\mu_i}:=x_i\restrict l_*$ for all  $i<N$ and let $r\in \mathbb S_\alpha$ be a condition such that 
$r\restrict \xi=u^{l_*}_{N-1}\restrict\xi$, $r\restrict\xi$ forces 
$r(\xi)$ to be $\Un\sset{u^{l_*}_i(\xi)}{i<N}$, and
$(r|\mu_i)\restrict\beta$  forces 
$r(\beta)=u^{l_*}_i(\beta)$ for all  $\beta>\xi$.
It follows from the above that $r|\mu_i\geq u^{l_*}_i$ and hence
$r|\mu_i$ forces $\tau\restrict l_*=x_i\restrict l_*=y_{\mu_i}$,
hence $q:=r$ and $l:=l_{*}$ are as required.

Now assume that $\card{F}>1$ and the statement holds for each set of cardinality smaller than $\card{F}$. Let $\xi,l_*,r$ be as above and note that by the construction we have
$$(r|\mu_i,n)\geq_{F\setminus\{\xi\}}(u^{l_*}_i,n) \geq_{F\setminus\{\xi\}}(p|\mu_i,n)$$ for all $i<N$,
and hence 
$(r,n)\geq_F (p,n)$.
Fix  $i<N$  and let $G$ be an $\bbS_{\xi+1}$-generic filter containing $(r|\mu_i)\restrict(\xi+1)$.
Work in $V[G]$. 
Then  
$$ r'_i:=\big((r|\mu_i)\restrict[\xi+1,\alpha)\big)^G  \in \bbS_{\xi+1,\alpha}
^G$$  
is $(F\sm\{\xi\},n)$-determined by Observation~\ref{obs:det2}(1) 
because 
$$(r'_i,n) \geq_{F\setminus\{\xi\}} \Big(\big(u^{l_*}_i\restrict[\xi+1,\alpha)\big)^{G},n\Big)$$
by the construction, and $\big(u^{l_*}_i\restrict[\xi+1,\alpha)\big)^{G}$ is $(F\sm\{\xi\},n)$-determined
by Observation~\ref{obs:cons3}.  Note that
\[
r_i'\forces \tau\in (\Cantor\cap V[\name{G}_{\xi+1,\alpha}])\sm\Un_{\xi<\beta<\alpha}(\Cantor\cap V[\name{G}_{\xi+1,\beta}]).
\]
because $r'_i\geq p\restrict[\xi+1,\alpha)^G$ and $p\restriction (\xi+1)\in G$. 
By the inductive assumption, there exist a condition $r''_i\in\bbS_{\xi+1,\alpha}^{G} $, a natural number $l_i>l_{*}$ and pairwise different elements $t_{\sigma'}\in 2^{l_i}$ for all maps $\sigma'\colon F\sm\{\xi\}\to 2^n$ consistent with $r_i''$,  such that $(r_i'',n)\geq_{F\sm\{\xi\}}(r_i',n)$ and $r_i''|\sigma'\forces \tau\restrict l_i = t_{\sigma'}$.
Let $\Sigma'_i$ be the set of all maps $\sigma'\colon F\sm\{\xi\}\to 2^n$ consistent with $r_i''$.

Now we work in $V$.
Let $\Sigma$ be the set of all maps $\sigma\colon F\to 2^n$ consistent with $p$.
Let $\underaccent{\tilde}{r}_i''$, $\underaccent{\tilde}{\Sigma'}_i$, 
$\underaccent{\tilde}{l}_i$ 
and $\underaccent{\tilde}{t}_{\sigma'}$ be $\bbS_{\xi+1}$-names for the condition $r_i''$, the set $\Sigma'_i$, natural number $l_i$ and finite sequences $t_{\sigma'}$, respectively.
Note that
\[
(r|\mu_i)\restrict(\xi+1)\forces \underaccent{\tilde}{\Sigma'}_i=\sset{\sigma\restrict F\sm\{\xi\}}{\sigma\in\Sigma\text{ and }\sigma(\xi)=s_i}
\]
by the second part of Observation~\ref{obs:cons3}.
By induction on $i<N$ pick a condition $r_i\in\mathbb S_{\xi+1}$ stronger than $(r|\mu_i)\restrict(\xi+1)$ and such that
$r_j\restrict\xi\geq r_i\restrict\xi$ for all $i<j<N$, $r_0\restrict\xi\geq r\restrict\xi$, and
which forces all the above properties, and also
decides all the  names mentioned in the previous sentences. More precisely, there exist $l_i>l_*$, $t^i_\sigma\in 2^{l_i}$ for all $\sigma\in \Sigma$ with $\sigma(\xi)=s_i$,
and $r_i''\in\mathbb S_{\xi+1,\alpha}$ such that $r_i$ 
forces that $\underaccent{\tilde}{r}_i''=r_i''$,
$\underaccent{\tilde}{l}_i=l_i$ and
$\underaccent{\tilde}{t}_{\sigma\restrict F\sm\{\xi\}}=t^i_{\sigma}$ for all maps $\sigma\in\Sigma$ with $\sigma(\xi)=s_i$. Thus  the elements $t^i_{\sigma}$ are pairwise different for all maps $\sigma\in\Sigma$
with $\sigma(\xi)=s_i$.

Let $w\in\mathbb S_\alpha$ be a condition such that 
$w\restrict\xi=r_{N-1}\restrict\xi$,
\[w(\xi)=\Un\sset{r_i(\xi)}{i<N},
\]
and for every ordinal number $\beta$ with $\xi<\beta<\alpha $ and natural number $i<N$, we have 
\[
(w|\mu_i)\restrict\beta \forces  w(\beta)=r_i''(\beta).
\]
It follows that  
$$w|\sigma\forces \tau\restrict l_i = t^i_\sigma$$
for all $\sigma\in\Sigma$
with $\sigma(\xi)=s_i$
because 
$$w|\sigma\geq r_i\cup r_i''|\big(\sigma\restrict(F\setminus\{\xi\})\big).$$

Let $\sigma,\nu\in\Sigma$ be different maps.
Assume that $\sigma(\xi)\neq\nu(\xi)$.
Then there are different natural numbers $i,j$ such that $\sigma(\xi)=s_i$ and $\nu(\xi)=s_j$.
Since $w|\sigma\geq w|\mu_i\geq u^{l_*}_i$ and 
$w|\nu\geq w|\mu_j\geq u^{l_*}_j$, we have
\[
w|\sigma\forces\tau\restrict l_i=t^i_\sigma,\quad\tau\restrict l_*\in \sset{x\restrict l_*}{x\in X_i}
\]
and
\[
w|\nu\forces\tau\restrict l_j=t^j_\nu\quad\tau\restrict l_*\in \sset{x\restrict l_*}{x\in X_j}.
\]
The sets $\sset{x\restrict l_*}{x\in X_i}$ and $\sset{x\restrict l_*}{x\in X_j}$ are disjoint, and thus $t^i_\sigma\restrict l_*\neq t^j_\nu\restrict l_*$.
Now assume that there is a natural number $i$ such that $\sigma(\xi)=\nu(\xi)=s_i$.
Then $\sigma\restrict F\sm\{\xi\}\neq \nu\restrict F\sm\{\xi\}$
and the condition $(w|\mu_i)\restrict(\xi+1)$ forces that $\underaccent{\tilde}{t}_{\sigma\restrict F\sm\{\xi\}}=t^i_\sigma$ and $\underaccent{\tilde}{t}_{\nu\restrict F\sm\{0\}}=t^i_\nu$
because $(w|\mu_i)\restrict(\xi+1)\geq r_i$, while  $t^i_\sigma\neq t^i_\nu$.
Summarizing, if $\nu(\xi)=\sigma(\xi)=s_i$, then $t^i_\sigma\neq t^i_\nu$, 
and if $\nu(\xi)=s_j\neq s_i=\sigma(\xi)$, then $t^i_\sigma\restrict l_*\neq t^j_\nu\restrict l_*$.

Finally, applying Observation~\ref{obs:det2}(7)
to $l=\max_{i<N}l_i$ and $w$ which is $(F,n)$-determined (because $(w,n)\geq_F(p,n)$ by the construction), we get a condition $q$  such that  $(q,n)\geq_F (w,n)$, and for every
$\sigma\in\Sigma$ a sequence $y_\sigma\in 2^l$  such that
$q|\sigma\Vdash\tau\restrict l=y_\sigma$.
Since $q|\sigma\geq w|\sigma,$
we have 
$y_\sigma\restrict l_i=t^i_\sigma$ for all $\sigma\in\Sigma$ with $\sigma(\xi)=s_i$.
It follows from the above that the $y_\sigma$'s are mutually different:
if $\nu(\xi)=\sigma(\xi)=s_i$, then $y_\sigma\restrict l_i =t^i_\sigma\neq t^i_\nu=y_\nu\restrict l_i$;  
and if $\nu(\xi)=s_j\neq s_i= \sigma(\xi)$, then $y_\sigma\restrict l_*=t^i_\sigma\restrict l_*\neq t^j_\nu\restrict l_*=y_\nu\restrict l_*$.
\epf

The following fact is reminiscent  of 
\cite[Lemma~2.3(i)]{BauLav}, and we use a rather similar approach 
to the proof, which we present for the sake of completeness.

\blem\label{lem:det}
Let $\alpha$  be an ordinal, $p\in\mathbb S_\alpha$, $n\in\w$  and $F\sub\alpha$  a nonempty finite set.
Then there are a natural number $k>n$ and an $(F,k)$-determined condition $q\in \mathbb S_\alpha$  such that $(q,k)\geq_F (p,n)$. 
\elem
\bpf
We proceed by induction on $|F|$. Suppose that $F=\{\beta\}$ for some 
$\beta<\alpha$ and pick $r\in\mathbb S_\beta$, $r\geq p\restrict\beta$
which decides $k>n$ and  $p(\beta)\cap 2^k$ (and hence also decides $p(\beta)\cap 2^n$), so that each element $s\in p(\beta)\cap 2^n$
has at least 2 extensions  in $p(\beta)\cap 2^k$.
Then $q:=r\cup p\restrict[\beta,\alpha)$ is as required.

Now assume that $|F|>1$ and let $\beta=\max(F)$. By the inductive assumption there exists $r\in\mathbb S_\beta$ and $k_0>n$ such that $r$ is
$(F\cap\beta,k_0)$-determined and $(r,k_0)\geq_{F\cap\beta} (p\restrict\beta,n)$.
Let $\Sigma$ be the family of all $\sigma:F\cap\beta\to 2^{k_0}$ consistent with $r$. Let $N=|\Sigma|$ and  write $\Sigma$ in the form $\{\sigma_i:i<N\}$.
By induction on $i<N$ let us construct a sequence
$$(r^0_{N-1},k_0)\geq_{F\cap\beta}\cdots \geq_{F\cap\beta}(r^0_{0},k_0)\geq_{F\cap\beta} (r^0_{-1},k_0),$$
where $r^0_{-1}=r$ and $r^0_i\in\mathbb S_\beta$ for all $i<N$,
as follows: Given  $r^0_{i-1}$ for some $i<N$, let $\mathbb S_\beta\ni u^0_i\geq r^0_{i-1}|\sigma_i$
be a condition such that there exists $T^0_i\sub  2^{l_i}$
for some $l_i\geq k_0$ such that 
$u^0_i$ forces $p(\beta)\cap 2^{l_i}=T^0_i$, and for every
$s\in T^0_i\restrict n$ there exists at least two $t\in T^0_i$ extending $s$.
Now, let $r^0_i\in\mathbb S_\beta$ be such that
$(r^0_i,k_0)\geq_{F\cap\beta}(r^0_{i-1},k_0)$ and
$r^0_i|\sigma_i=u^0_i$, by Observation~\ref{obs:det2}(5).

Let $k=\max_{i<N}l_i$ and set $r^1_{-1}=r^0_{N-1}$.
By induction on $i<N$ let us construct a sequence
$$(r^1_{n-1},k_0)\geq_{F\cap\beta}\cdots \geq_{F\cap\beta}(r^1_{0},k_0)\geq_{F\cap\beta} (r^1_{-1},k_0),$$
where  $r^1_i\in\mathbb S_\beta$ for all $i<N$,
as follows: Given  $r^1_{i-1}$ for some $i<N$, let $\mathbb S_\beta\ni u^1_i\geq r^1_{i-1}|\sigma_i$
be a condition such that there exists $T^1_i\sub  2^{k}$
for which 
$u^1_i$ forces $p(\beta)\cap 2^{k}=T^1_i$. 
Now, let $r^1_i\in\mathbb S_\beta$ be such that
$(r^1_i,k_0)\geq_{F\cap\beta}(r^1_{i-1},k_0)$ and
$r^1_i|\sigma_i=u^1_i$.
Note that $T^0_i=\{t\restrict l_i:t\in T^1_i\}$.

Set $r^2_{-1}=r^1_{N-1}$.
By induction on $i<N$ let us construct a sequence
$$(r^2_{n-1},k_0)\geq_{F\cap\beta}\cdots \geq_{F\cap\beta}(r^2_{0},k_0)\geq_{F\cap\beta} (r^2_{-1},k_0),$$
where  $r^2_i\in\mathbb S_\beta$ for all $i<N$,
as follows: Given  $r^2_{i-1}$ for some $i<N$, let $\mathbb S_\beta\ni u^2_i\geq r^2_{i-1}|\sigma_i$
be a condition such that 
for every $\gamma\in F\cap\beta$ there exists $\nu_i(\gamma)\in 2^k$
such that $\sigma_i(\gamma)=\nu_i(\gamma)\restrict k_0$
and $u^2_i\restrict\gamma$ forces that $\nu_i(\gamma)$ is an initial segment of the stem of $u^2_i(\gamma)$. Such an $u^2_i$ can be constructed recursively over $\gamma\in F\cap\beta$, moving from the bigger to smaller elements. 
Now, let $r^2_i\in\mathbb S_\beta$ be such that
$(r^2_i,k_0)\geq_{F\cap\beta}(r^2_{i-1},k_0)$ and
$r^1_i|\sigma_i=u^2_i$.

We claim that $q=r^2_{N-1}\cup p\restrict[\beta,\alpha)$ and $k$
are as required. Indeed, we have that 
$(q\restrict\beta,k)\geq_{F\cap\beta} (p\restrict\beta,n)$
because $q\restrict\beta=r^2_{N-1}$, $k\geq k_0$, and 
$(r^2_{N-1},k_0)\geq_{F\cap\beta} (r,k_0)\geq_{F\cap\beta} (p\restrict\beta,n)$.
Moreover, since $r^2_{N-1}|\sigma_i\geq u^0_i$, we have that
$r^2_{N-1}|\sigma_i$ decides $p(\beta)\cap 2^k$ as $T^1_i$,
which has the property that any $s\in T^1_i\restrict n$ has at least two extensions
in $T^1_i$ (because $T^0_i=T^1_i\restrict l_i$ has this property).
It follows that $r^2_{N-1}|\sigma_i$ forces $q(\beta)=p(\beta) $ and $(p(\beta),k)\geq (p(\beta),n)$. Since
$\{r^2_{N-1}|\sigma_i:i<N\}$ is dense above $r^2_{N-1}$,
we conclude that $r^2_{N-1}=q\restrict\beta$ forces 
$(q(\beta),k)\geq (p(\beta),n)$, and therefore $(q,k)\geq (p,n)$.

Finally, by the construction of $r^2_{N-1}$ we have that 
$q$ is $(F,k)$-determined, with 
$$\big\{\nu_i\cup\{\la\beta,t\ra\}\: :\: i<N,t\in T^1_i\big\}$$
being the family of those $\nu:F\to 2^k$ which are consistent with $q$.
\epf

\blem[{Miller~\cite[Lemma~6]{mill}}]\label{lem:Mill}
Let $\alpha$ be an ordinal number, $p_0\in\mathbb S_\alpha$,   and $\tau$  an $\mathbb S_\alpha$-name for a real   such that
\[
p_0\forces \tau\in \bigl(\Cantor\cap V[\name{G}_{\alpha}]\bigr)\sm \Un_{\beta<\alpha}\bigl(\Cantor\cap V[\name{G}_\beta]\bigr).
\]
Then there exist a condition $p\geq p_0$, an increasing sequence $\seq{F_n}{n\in\w}$ of finite subsets of $\alpha$, increasing sequences of natural numbers $\seq{k_n}{n\in\w}$, $\seq{l_n}{n\in\w}$, and elements $y_{\sigma}\in 2^{l_n}$ for all maps $\sigma\colon F_n\to 2^{k_n}$ consistent with $p$, with the following properties:
\be
\item $\Un_{n\in\w}F_n=\supp(p)$,
\item $p$ is $(F_n,k_n)$-determined,
\item $(p, k_{n+1}) \geq_{F_n} (p, k_n)$,
\item $p|\sigma\forces \tau\restrict l_n= y_{\sigma}$
for all $\sigma\in (2^{k_n})^{F_n}$ consistent with $p$, and 
\item the maps $y_{\sigma}$, where $\sigma$ is as above, are mutually different.
\ee 
\elem

\bpf
The choice of the $F_n$'s is standard and thus will not be 
specified, except that we set $F_0=\{0\}$. Set also $k_0=0$. 
Trivially, $p_0$ is $(F_0,k_0)$-determined  since the unique  map $\{\la 0,\emptyset\ra\} $ in $(2^{k_0})^{F_0}$   is consistent with $p_0$.
By Lemma~\ref{lem:stepmill}, there are a condition $q_0\in\mathbb S_\alpha$, a natural number $l_0$ and pairwise different elements $y_{\sigma}\in 2^{l_0}$ for all maps\footnote{As we noted above, there is just one   map like that.} $\sigma\colon F_0\to 2^{k_0}$ consistent with $p_0$ such that $(q_0,k_0)\geq_{F_0}(p_0,k_0)$ and $q_0|\sigma\forces \tau\restrict l_0=y_{\sigma}$.

Let $k_1>k_0$ be a natural number and $p_1\in\mathbb S_\alpha$  a condition from Lemma~\ref{lem:det}, applied to the set $F_1$ and the condition $q_0$.
Then $p_1$ is $(F_1,k_1)$-determined and $(p_1,k_1)\geq_{F_0}(p_0,k_0)$.

Fix a natural number $n\geq 1$ and assume that a set $F_n$, natural number $k_n$, and an $(F_n,k_n)$-determined condition $p_n\in\mathbb S_\alpha$
with $(p_n,k_n)\geq_{F_{n-1}}(q_{n-1},k_{n-1})$ 
have already been defined.
By Lemma~\ref{lem:stepmill}, there are a condition $q_n\in\mathbb S_\alpha$, a natural number $l_n>l_{n-1}$ and pairwise different elements $y_{\sigma}\in 2^{l_n}$ for all maps $\sigma\colon F_n\to 2^{k_n}$ consistent with $p_n$ such that $(q_n,k_n)\geq_{F_n}(p_n,k_n)$ and $q_n|\sigma\forces \tau\restrict l_n=y_{\sigma}$.

Let $p$ be the fusion of the sequence $\seq{(p_n,k_n,F_n)}{n\in\w}$~\cite[Lemma~1.2]{BauLav} and note that it is as required.
\epf

Let $F$ be a subset of $H$, $n\leq m$ be natural numbers, and $\nu\colon F\to 2^n$, $\sigma\colon H\to 2^m$ be maps.
Following our convention at the beginning of Section~\ref{sec2}, the map $\sigma$ is an \emph{extension}
of $\nu$ (we denote this by $\nu\prec\sigma$) if $\nu(\beta)=\sigma(\beta)\restrict n$ for all $\beta\in F$. The next fact is standard.

\bobs\label{obs:ext_map}
In the notation above, if $F,H\sub \alpha$, $p\in\mathbb S_\alpha$,
$\nu,\sigma$ are consistent with $p$ and $\nu\prec\sigma$, then 
$p|\sigma\geq p|\nu$.
\eobs

\blem\label{lem:prompt}
Let $p\in\mathbb S_\alpha$ be an $(F,n)$-determined condition and $\Sigma_n$  be the set of all maps $\nu\colon F\to 2^n$ consistent with $p$. 
Then for every $G\sub F$ and $k\leq n$, if $p$
is $(G,k)$-determined and 
$\mu:G\to 2^k$ is consistent with $p$, then there exists
$\nu\in\Sigma_n$ extending $\mu$. 

Moreover, if $H\sub\alpha$ is a finite set with $F\sub H$,   $\beta\in F$, $m>n$ is a natural number,  $p$ is also $(H,m)$-determined,
$\Sigma_m$ is the set of all maps $\sigma\colon H\to 2^m$ consistent with $p$, and
$(p,m)\geq_F (p,n)$, then the following assertions hold.
\begin{itemize}
\item[(1)] 
 For every $\nu\in \Sigma_n$ and $\sigma\in\Sigma_m$ 
such that $\sigma\restrict (H\cap \beta)$ extends $\nu\restrict (F\cap\beta)$,
there are $\sigma_1,\sigma_2\in\Sigma_m$ extending $\nu$ with $\sigma_1\restrict (H\cap \beta)=\sigma_2\restrict (H\cap \beta)=\sigma\restrict (H\cap\beta)$, and  such that $\sigma_1(\beta), \sigma_2(\beta)$
are distinct extensions of $\nu(\beta)$.
\item[$(2)$] For every $\nu\in \Sigma_n$ and $\rho:H\cap\beta\to 2^m$ consistent with $p$, if $\rho$ extends $\nu\restrict (F\cap\beta)$,
then there are $\sigma_1,\sigma_2\in\Sigma_m$ extending $\nu$ with $\sigma_1\restrict (H\cap \beta)=\sigma_2\restrict (H\cap \beta)=\rho$ and  such that $\sigma_1(\beta), \sigma_2(\beta)$
are distinct extensions of $\nu(\beta)$.
\item[$(3)$]
For every $C\subset\Sigma_n$, coherent $e_0:C\to \Sigma_m$, 
and $\nu_*\in\Sigma_n\setminus C$, 
there exist two different coherent maps $e,e'\colon (C\cup\{\nu_*\}) \to\Sigma_m$ such that
$e\restrict C=e'\restrict C=e_0$.
\end{itemize}
\elem
\bpf
We start with proving the first part.
 Proceed by induction on $|G|$. If $G=\emptyset$, then there is nothing to prove.
Let $\gamma:=\max G$ and assume that the statement holds for $G':=G\setminus\{\gamma\}$.
Fix a map $\mu\colon G\to 2^k$ consistent with $p$ and let $\mu':=\mu\restrict G'$.
By Observation~\ref{obs:det2}(4), $p$ is both $(F\cap\gamma,n)$- and $(G\cap\gamma,k)$-determined, and hence there exists 
 $\nu':(F\cap\gamma)\to  2^n$ consistent with $p$ such that  $\nu'\succ\mu'$.
Since $(p\restrict\gamma)|\mu'\forces \mu(\gamma)\in p(\gamma),$
the condition $(p\restrict\gamma)|\nu'$ also forces this because
$(p\restrict\gamma)|\nu'\geq (p\restrict\gamma)|\mu'$ by Observation~\ref{obs:ext_map}.
Strengthening the latter condition to some  $r$ if necessary, we may find $t\in 2^n$ extending  $\mu(\gamma)$
such that $r\Vdash t\in p(\gamma)$.
Since $p$ is $(F,n)$-determined, we conclude that
$(p\restrict\gamma)|\nu'\forces t\in p(\gamma).$
It follows that $\nu'\cup\{(\gamma,t)\}\in (2^n)^F$ is consistent with $p$ and extends $\mu$.

The second part of the lemma is rather straightforward.
Item  $(2)$ is an equivalent reformulation of $(1)$, by Observation~\ref{obs:det2}(2).  We shall present the proof of $(3)$, because we find it the least obvious one.

 Proceed by induction on $|F|$. 
 Assume that $F=\{\beta\}$. 
Let $\nu \in C$ and $\sigma := e_0(\nu)$. 
Choose  some distinct $\mu,\mu'\in\Sigma_m$ such that $\mu,\mu'\succ\nu_*$ and $\mu \restriction (H \cap \beta) = \sigma \restriction (H \cap \beta) = \mu' \restriction (H \cap \beta)$, which is possible by (1). 
Then the maps $e=e_0\cup\{\la\nu_*,\mu\ra\}$ and
$e'=e_0\cup\{\la\nu_*,\mu'\ra\}$ are easily seen to be as required. 

Now assume that $|F|>1$ and $(3)$ holds for any finite subset of the support of $p$ of cardinality smaller than $\card{F}$.
Set $\beta=\max F$,
   $C^-=\{\nu\restrict(F\cap\beta):\nu\in C\}$, 
 $e^-_0(\nu\restrict(F\cap\beta))= e_0(\nu)\restrict(H\cap\beta)$ for all $\nu\in C$,  
   and $\nu^-_*=\nu_*\restrict(F\cap\beta)$. Let also
   $\Sigma_m^-$ be the family of all maps $\sigma:H\cap\beta\to 2^m$ consistent with $p$.
 By the inductive assumption  applied to the objects defined above 
 we can get a coherent 
 $e^-:C^-\cup\{\nu^-_*\}\to \Sigma^-_m $
 such that $e^-\restrict C^-= e^-_0$.
 (Actually, we could get even two different such $e^-$  if $\nu^-_*\not\in C^-$,
 but thus irrelevant here.)

By $(2)$ applied to $\nu_*$ and $\rho=e^-(\nu^-_*)$ we
can find distinct $\sigma,\sigma'\in\Sigma_m$ such that 
$$\sigma\restrict (H\cap\beta)=\rho=\sigma'\restrict (H\cap\beta).$$
It suffices to show that $e=e_0\cup\{\la\nu_*,\sigma\ra\}$
and $e'=e_0\cup\{\la\nu_*,\sigma'\ra\}$ are both coherent. We shall check this for
$e$, the case of $e'$ is analogous. 
Pick $\nu\in C$  and let
$\gamma\in F$ be the minimal element with $\nu(\gamma)\neq\nu_*(\gamma)$.
Such an ordinal $\gamma$ must exist because $\nu\neq\nu_*$ as $\nu\in C\not\ni\nu_*$.
Fix $\beta\geq \gamma$.
In the assertions from the cases below we use the fact that
$$e^-(\nu_*\restrict (F\cap\beta))=e^-(\nu^-_*)=\rho =\sigma\restrict (H\cap\beta)$$
and
$e(\nu_*)=\sigma$. 

If $\gamma<\beta$, then
\begin{multline*}
 e(\nu)\restrict (H\cap\gamma)= \big(e(\nu)\restrict (H\cap\beta)\big) \restrict (H\cap\gamma)= \big(e^-(\nu\restrict (F\cap\beta)\big) \restrict (H\cap\gamma) = \\ 
 = \big(e^-(\nu_*\restrict (F\cap\beta)\big) \restrict (H\cap\gamma)=
 \big(e(\nu_*)\restrict (H\cap\beta)\big) \restrict (H\cap\gamma)=
 e(\nu_*)\restrict(H\cap\gamma).
\end{multline*}
    If $\gamma=\beta$, then
\[
 e(\nu)\restrict (H\cap\beta)= 
 e^-(\nu\restrict (F\cap\beta)) 
 = e^-(\nu_*\restrict (F\cap\beta)) =e(\nu_*)\restrict (H\cap\beta).\qedhere
\]
\epf

\blem \label{clm:Sigma}
We use notation  from the formulation of Lemma~\ref{lem:Mill}. 
Set $S=\mathrm{supp}(p)$, $\Sigma_n=\{\nu\in (2^{k_n})^{F_n}:\nu$ is consistent with $p\}$, and 
suppose that 
$\seq{\la i_n,j_n,C_n \ra}{n\in\w}$ is  a sequence such that
items $(1)$--$(4)$ of Lemma~\ref{clm:easy} are satisfied.
Then there exists $q\geq p$ such that 
$\{p|\sigma:\sigma\in C_n\}$ is predense above $q$ for all $n\in\w$.
\elem
\bpf
For each natural number $n$,
let $q_n\geq p$ be such that for every 
$\beta<\alpha$ and $\nu\in C_n$ we have 
\[
(q_n\restrict\beta)|(\nu\restrict (F_{i_n}\cap \beta))\forces_\beta q_n(\beta)=
\bigcup \smallmedset{(p|\nu')(\beta)}{\nu'\in C_n,\nu'\restrict (F_{i_n}\cap \beta)= \nu\restrict (F_{i_n}\cap \beta)}.
\]
The correctness of this definition formally
requires to prove by recursion over $\beta\leq\alpha$, along with defining $q_n$, that 
each $\nu\restrict (F_{i_n}\cap \beta)$ is consistent with $q_n\restrict\beta$, 
where $\nu\in C_n$,
and the set
\[ \smallmedset{(q_n\restrict\beta)|(\nu\restrict (F_{i_n}\cap \beta))}{\nu\in C_n}
\]
is a maximal antichain above $q_n\restrict\beta$, but this is rather easy and standard.
As a result, the set $ \sset{q_n|\nu}{\nu\in C_n}$
is a maximal antichain above $q_n$.

It remains to show that 
\[
(q_{n+1},k_{i_{n+1}})\geq_{F_{i_n}}(q_n,k_{i_n})
\]
for all $n\in\w$, and then let $q$  be the fusion of the $q_n$'s. 
Suppose that for some 
$\beta\in F_{i_n}$ we have already shown that
\[
(q_{n+1}\restrict\beta,k_{i_{n+1}})\geq_{F_{i_n}\cap\beta}(q_n\restrict\beta,k_{i_n}),
\]
and we will prove that
\[
q_{n+1}\restrict\beta\Vdash (q_{n+1}(\beta), k_{i_{n+1}})\geq (q_n(\beta), k_{i_n}),
\]
This boils down to proving that for every $\sigma\in C_{n+1}$,
if $\nu(\gamma):=\sigma(\gamma)\restrict k_n$ for all $\gamma\in F_{i_n}$,
then
\begin{multline*}
(q_{n+1}\restrict\beta)|(\sigma\restrict (F_{i_{n+1}}\cap \beta))\forces_\beta \\
\Bigl(\Un \smallmedset{(p|\sigma')(\beta)}{\sigma'\in C_{n+1},\sigma'\restrict (F_{i_{n+1}}\cap \beta)= \sigma\restrict (F_{i_{n+1}}\cap \beta)}, k_{n+1}\Bigr)\geq \label{long_eq}  \\
\geq
\Bigl(\Un\smallmedset{(p|\nu')(\beta)}{\nu'\in C_n,\nu'\restrict (F_{i_n}\cap \beta)= \nu\restrict (F_{i_n}\cap \beta)},k_n\Bigr).
\end{multline*}
Fix $s\in 2^{k_n}$ such that there exists $\nu'\in C_n$
with $\nu'(\beta)=s$  and $\nu'\restrict (F_{i_n}\cap\beta)=\nu\restrict (F_{i_n}\cap\beta)$.  By Lemma~\ref{clm:easy}(4), we have $e(\nu)\restrict (F_{j_n}\cap\beta)=
e(\nu')\restrict (F_{j_n}\cap\beta)$.
By Lemma~\ref{lem:prompt}(1), there are $\sigma_1,\sigma_2\in\Sigma_{i_{n+1}}$, both
extending $e(\nu')$ (and hence $\sigma_1,\sigma_2\in C_{n+1}$) such that  $\sigma_1(\beta)\neq\sigma_2(\beta)$ and 
\[
\sigma_1\restrict(F_{i_{n+1}}\cap\beta)=\sigma_2\restrict(F_{i_{n+1}}\cap\beta)=\sigma\restrict(F_{i_{n+1}}\cap\beta).
\]
It follows that 
\[
(q_{n+1}\restrict\beta)|(\sigma\restrict(F_{i_{n+1}}\cap\beta))\forces_\beta\sigma_1(\beta),\sigma_2(\beta)\in q_{n+1}(\beta)\cap 2^{k_{i_{n+1}}}
\]
and $\sigma_1(\beta),\sigma_2(\beta)$ are distinct extensions 
of $e(\nu')(\beta),$ which in its turn extends $\nu'(\beta)=s$.
\epf



\section{Totally imperfect Menger sets and Sacks forcing } \label{menger_size}

By $V$ we mean a ground model of ZFC and $G_{\w_2}$ is an $\mathbb S_{\w_2}$-generic filter over $V$.

\bthm\label{thm:main_ti}\mbox{}
In $V[G_{\w_2}]$, every totally imperfect Menger set
$X\sub 2^\w$ has size at most $\w_1$. 
 \ethm

For the proof of Theorem~\ref{thm:main_ti}, we need the following auxiliary result, whose proof is rather standard (see, e.g., 
\cite[Lemma~5.10]{BlaShe87} for a similar argument) and is left to the reader.

\blem\label{lem:LS_ti}
In $V[G_{\w_2}]$, let $X\sub 2^\w$. Then there exists $\alpha<\w_2$ of cofinality $\w_1$ such that
\be 
\item $X\cap V[G_\alpha]\in V[G_\alpha]$ and if $K,K'\sub 2^\w$ are closed 
crowded sets and coded in $V[G_\alpha]$, and  $K'\sub K\setminus (X\cap V[G_\alpha])$,
then $K'\sub K\setminus X$.
\ee
Moreover, if $X$ is a 
totally imperfect Menger set in $V[G_{\w_2}]$, then 
\begin{itemize}
\item[(2)] $X\cap V[G_\alpha]$ is a totally imperfect Menger set in $V[G_\alpha]$.
\end{itemize}
\elem

In what follows we use the same notation for a Borel (typically closed) subset of $ 2^\w$
in the ground model as well as for its reinterpretation in the forcing extensions we consider, it will be always clear in which set-theoretic universe we work.

\bpf[{Proof of Theorem~\ref{thm:main_ti}}]
Let $\alpha$ be such as in Lemma~\ref{lem:LS_ti}. Working in $V[G_{\w_2}]$,
we claim that $X\sub  V[G_\alpha]$.
Since in $V[G_\alpha]$ the remainder 
$\mathbb S_{\alpha,\w_2}$ is order-isomorphic to $\mathbb S_{\w_2}$, there is no loss of generality in  
assuming that $\alpha=0$, i.e., that $V=V[G_\alpha]$. 

Let us pick $z\in X\setminus V$ and let $\gamma$ be the minimal ordinal with 
$z\in V[G_\gamma]$. Next, we work in $V$.
Let $p_0\in G_{\gamma}$ and $\tau\in V^{\mathbb S_\gamma}$ be 
such that $\tau^{G_\gamma}=z$ and
$$p_0\Vdash_\gamma\tau \in (2^\w \cap V[\dot{G}_\gamma])\setminus\bigcup_{\beta<\gamma}V[\dot{G}_\beta].$$ 
We shall find $q\geq p_0$, $q \in \mathbb{S}_\gamma$ such that $q\Vdash_{\w_2}\tau\not\in X$. This would accomplish the proof: The genericity
of $G_{\gamma}$ implies that there is $q$ as above which lies in $G_{\gamma} \subseteq G_{\w_2}$, 
which would yield $z=\tau^{G_\gamma}=\tau^{G_{\w_2}}\not\in X$. 

Take $p$ and $F_n,k_n,l_n, y_{\sigma_n}$ from Lemma~\ref{lem:Mill}, applied to $p_0$ and  $\tau$. Note also that Lemma~\ref{lem:prompt}(3) ensures that 
$S=\mathrm{supp}(p)$ and the sequences $\seq{k_n}{n\in\w}, \seq{F_n}{n\in\w}, \seq{\Sigma_n}{n\in\w}$
satisfy $(e_f)$ from the first paragraph of Section~\ref{sec2}.
Let $K$ and $[\sigma]$ for $\sigma\in\bigcup_{n\in\w}\Sigma_n$
be defined in the same way as in the first paragraph of Section~\ref{sec2}.

Fix an element $x\in K$ and let $\sigma_n\in\Sigma_n$ be  such that $\{x\}=\bigcap_{n\in\w}[\sigma_n]$.
Fix a natural number $n$.
We have $\sigma_n(\beta)\sub\sigma_{n+1}(\beta)$ for all $\beta\in F_n$, i.e., $\sigma_n \prec \sigma_{n+1}$.
Then $p|\sigma_{n+1}\geq p|\sigma_n$ and it follows from Lemma~\ref{lem:Mill}(4) that
\[
p|\sigma_{n+1}\forces \tau\restriction l_n=y_{\sigma_n}\text{ and }p|\sigma_{n+1}\forces \tau \restriction l_{n+1}=y_{\sigma_{n+1}},
\]
which gives  $y_{\sigma_n}\sub y_{\sigma_{n+1}}$.
Thus, the map $h\colon K\to \Cantor$ such that
\[
h(x):=\Un_{n\in\w}y_{\sigma_n},
\]
for all $x\in K$, is well defined.
By Lemma~\ref{lem:Mill}(5), the map $h$ is a continuous injection. 
Consequently, the map $h\colon K\to h[K]$ is a homeomorphism, and hence $h[K]$ is perfect.

Fix a natural number $n$.
By Lemma~\ref{lem:Mill}(2) and Observation~\ref{obs:det2}(3), the set $\sset{p|\sigma}{\sigma\in\Sigma_n}$  is a maximal antichain above $p$.
Applying Lemma~\ref{lem:Mill}(4), we have that 
\[
p\forces \tau\restrict l_n\in \sset{y_{\sigma}}{\sigma\in\Sigma_n}.
\]
It follows from the above and from the definition of the function $h$  that
$p\forces \tau\in h[K].$ 

By our assumption on $\alpha=0$, the set $h[K]\cap X\cap V$ is an element of $V$ and it is totally imperfect and Menger in $V$.
Let $\seq{\la i_n,j_n,C_n\ra}{n\in\w}$ be a sequence from Lemma~\ref{clm:easy}, applied to $S:=\mathrm{supp}(p)$ and $h^{-1}[X\cap V]\sub K$.
Lemma~\ref{clm:easy}(5) yields
\[
K':=\bigcap_{n\in\w}\Un_{\nu\in C_n}[\nu]\sub K\sm h\inv[X\cap V],
\]
and therefore, $h[K']\sub h[K]\setminus (X\cap V)$.
Since $K,K',h$ are all coded in $V$, we conclude that
$h[K']\sub h[K]\setminus X$ holds in $V[G_{\w_2}]$ by Lemma~\ref{lem:LS_ti}(2).

Let $q\geq p$ be a condition given by Lemma~\ref{clm:Sigma}.
Since the sets $\{p|\nu:\nu\in C_n\}$ are predense above $q$ for all $n\in\w$, we have $q\Vdash\tau \restriction l_{i_n} \in \sset{y_{\nu}}{\nu\in C_n}$, and thus
$q\Vdash\tau\in h[K']$.
We conclude that 
$q\Vdash \tau\not\in X$.
\epf


\section{A modification of the Menger game and consonant spaces} \label{consonant}

Let $X\sub\Cantor$.
We introduce a modification of the Menger game played on $X$, which we call 
\emph{grouped Menger game} played on $X$.

Round 0: 
\Alice{} selects a natural number $l_0>0$, and then the players 
play the usual Menger game $l_0$ subrounds, thus constructing a partial play
\[
(l_0,\cU_0, \cF_0,\dotsc,\cU_{l_0-1},\cF_{l_0-1}),
\]
where $\cF_i$ is a finite subfamily of $\cU_i$ for all natural numbers $i<l_0$.

Round 1:
\Alice{} selects a natural number $l_1>0$, and then the players 
play the usual Menger game additional $l_1$ subrounds, thus constructing a partial play
\[
(l_0,\mathcal U_0,\ldots,\mathcal U_{l_0-1},\mathcal F_{l_0-1};\: l_1, \mathcal U_{l_0},\mathcal F_{l_0},\ldots,\mathcal U_{l_0+l_1-1},\mathcal F_{l_0+l_1-1}),
\]
where $\cF_i$ is a finite subfamily of $\cU_i$ for all natural numbers $i<l_0+l_1$.

Fix a natural number $n>0$ and assume that natural numbers $l_0,\dotsc,l_{n-1}>0$ and $n-1$ rounds of the game have been defined.

Round $n$: 
\Alice{} selects a natural number $l_n>0$, and then the players 
play the usual Menger game additional $l_n$ subrounds,
thus constructing a partial play
\begin{multline} 
(l_0,\mathcal U_0,\cF_0,\dotsc,\mathcal U_{l_0-1},\mathcal F_{l_0-1};\  l_1, \mathcal U_{l_0},\mathcal F_{l_0},\dotsc,\mathcal U_{l_0+l_1-1},\mathcal F_{l_0+l_1-1};\dotsc \nonumber ;\\
l_n, \mathcal U_{l_0+l_1+\cdots+l_{n-1}},\mathcal F_{l_0+l_1+\cdots+l_{n-1}},\ldots,  \mathcal U_{l_0+l_1+\cdots+l_{n}-1},\mathcal F_{l_0+l_1+\cdots+l_{n}-1}), \label{play_look}    
\end{multline}
where $\cF_i$ is a finite subfamily of $\cU_i$ for all natural numbers $i<l_0+\dotsb+l_{n}$.

Let $L_0:=0$
and $L_{n+1}:=l_0+l_1+\dotsb+l_n$ for all natural numbers $n$.
\Bob{} wins the game if 
\[
X=\bigcup_{n\in\w}\bigcap_{i\in [L_n,L_{n+1})}\Un\mathcal F_n,
\]
and \Alice{} wins otherwise. 

\brem\label{rem:gMincr}
In the grouped Menger game played on a set $X\sub\Cantor$, there is no loss of generality if we assume that covers given by \Alice{} in each step are countable and increasing and the families $\cF_n$ chosen by \Bob{} are singletons.
\erem

We are interested in sets $X\sub 2^\w$ for which \Alice{}
has no winning strategy in the grouped Menger game played on $X$.
They include two important 
classes of subspaces of the Cantor space: Hurewicz spaces and those ones whose complement is consonant.

The following observation is an immediate consequence of Theorem~\ref{thm:hur_hur}.

\bobs\label{obs:hur}
If $X\sub 2^\w$ is Hurewicz, then 
\Alice{} has no winning strategy in the grouped Menger game played on $X$.
\eobs

\bprp\label{obs_obs}
If $Y\sub 2^\w$ is consonant, then 
\Alice{} has no winning strategy in 
the grouped Menger game played on $2^\w\sm Y$.
\eprp
\bpf
Let $\sigma$ be a strategy for \Alice{} in the grouped Menger game played on $X:=\Cantor\sm Y$.
By Remark~\ref{rem:gMincr}, we may assume that each family given by \Alice{} according to the strategy $\sigma$ is countable and increasing and \Bob{} chooses one set from the families given by \Alice{}. 
We shall define a strategy $\S$ for \Alice{} in the game $\GKO$ played on $X$ such that each play in $\GKO$ played according to $\S$ and lost by \Alice{}, gives rise
to a play in the grouped Menger game on $X$ in which \Alice{} uses 
$\sigma$ and loses. 

Suppose that $\sigma$ instructs
\Alice{} to start round $0$ by selecting a natural number $l_0>0$.
Let $L_0:=0$ and $L_1:=l_0$.
Then $\S$ instructs \Alice{} to play the family of all sets 
\[
\bigcap_{i\in[L_0,L_1)}U_i,
\]   
where
\[
(l_0,\mathcal U_0,U_0,\ldots,\mathcal U_{l_0-1},U_{l_0-1})
\]
is a play, where \Alice{} uses the strategy $\sigma$.
Note that the family of all 
the intersections as above is indeed an open $k$-cover of $X$.
Suppose that \Bob{} replies in  $\GKO$
by selecting $\bigcap_{i\in[L_0,L_1)}U_i$ for some sequence $(\mathcal U_0,U_0,\dotsc,\mathcal U_{l_0-1},U_{l_0-1})$ as above.
The strategy $\sigma$ instructs
\Alice{} to proceed in  round $1$ by selecting a natural number $l_1>0$.
Let $L_2:=l_0+l_1$.
Then $\S$ instructs \Alice{} to play the family of all sets
\[
\bigcap_{i\in[L_1,L_2)}U_i,
\]  
where
\[
(l_0,\mathcal U_0,U_0,\dotsc,\mathcal U_{l_0-1},U_{l_0-1};\: l_1, \mathcal U_{l_0},U_{l_0},\ldots,\mathcal U_{l_0+l_1-1},U_{l_0+l_1-1})
\] 
is a play in which \Alice{} uses $\sigma$, an open $k$-cover of $X$.
Suppose that \Bob{} replies 
in  $\GKO$
by selecting $\bigcap_{i\in[L_1,L_2)}U_i$
for some sequence $(\mathcal U_0,\mathcal F_0,\dotsc,\mathcal U_{l_0-1},U_{l_0-1};\: \mathcal U_{l_0},U_{l_0},\ldots,\allowbreak\mathcal U_{l_0+l_1-1},U_{l_0+l_1-1})$ as above.
\smallskip

In general, let $\sigma$ instruct
 \Alice{} to start round $n$ by selecting a natural number $ l_n>0$. 
Let $L_{n+1}:=l_0+l_1+\dotsb +l_n$.
Then the next move of \Alice{}  in  $\GKO$
according to $\S$ is, by the definition, the family of all sets
\[
\bigcap_{i\in[L_n,L_{n+1})} U_i,
\]
where
\begin{multline*}
 ( l_0,\mathcal U_0,U_0,\ldots,\mathcal U_{ l_0-1},U_{ l_0-1};\   l_1, \mathcal U_{ l_0},U_{ l_0},\ldots,\mathcal U_{ l_0+ l_1-1},U_{ l_0+ l_1-1};\ldots
 ;\\
 l_n, \mathcal U_{ l_0+ l_1+\cdots+ l_{n-1}},U_{ l_0+ l_1+\dotsb+ l_{n-1}},\ldots,  \mathcal U_{ l_0+ l_1+\dotsb+ l_{n-1}+ l_n-1},U_{ l_0+ l_1+\dotsb+ l_{n-1}+ l_n -1})
\end{multline*}
is a play in which \Alice{} uses $\sigma$, an open $k$-cover of $X$.

Since the strategy 
$\S$ is not winning, there is a play in  $\GKO$
in which \Alice{} uses $\S$ and loses, which gives rise to an infinite play
\begin{multline*} 
 ( l_0,\mathcal U_0,U_0,\ldots,\mathcal U_{ l_0-1},U_{ l_0-1};\   l_1, \mathcal U_{ l_0},U_{ l_0},\ldots,\mathcal U_{ l_0+ l_1-1},U_{ l_0+ l_1-1}; \ldots;\\
 l_n, \mathcal U_{ l_0+ l_1+\cdots+ l_{n-1}},U_{ l_0+ l_1+\cdots+ l_{n-1}},\ldots,  \mathcal U_{ l_0+ l_1+\cdots+ l_{n-1}+ l_n-1},U_{ l_0+ l_1+\cdots+ l_{n-1}+ l_n -1};  l_{n+1},\ldots)
\end{multline*}
 in which \Alice{} uses $\sigma$ and 
\[
X=\bigcup_{n\in\w}\bigcap_{i\in[L_n,L_{n+1})} U_i .
\]
This means that the strategy $\sigma$ is not winning as well.
\epf

\subsection{Menger game versus the grouped Menger game}

Let $\mathcal{GM}$ be the class of all subspaces $X$ of $2^\w$
such that \Alice{} has no winning strategy in the grouped Menger game played on 
$X$. Obviously, $\mathcal{GM}$ is contained in the class of all Menger subspaces of $2^\w$.
As we established in Section~\ref{consonant}, $\mathcal{GM} $ includes 
Hurewicz subspaces and subspaces with consonant complement, and hence also 
all Rothberger subspaces of $2^\w$ (see the discussion at the end of Section 1 in the work of Jordan~\cite{Jor20}).
Our next result gives a consistent example of a Menger space which does 
not belong to $\mathcal{GM}$.

\bprp \label{ultraf}
The class $\mathcal{GM}$ contains no ultrafilters.
\eprp
\bpf
Given an ultrafilter $X$ on $\w$, we shall describe a winning strategy  $\sigma$ for \Alice{} in the grouped Menger game played on $X$.
For natural numbers $n<k$, let 
\[
U_{[n,k)}:=\sset{a\sub \w}{a\cap [n,k)\neq\emptyset}.
\]
Then the families $\mathcal U_n:=\sset{U_{[n,k)}}{k> n}$ are increasing open covers of $X$ for all $n\in\w$.
Playing according to the strategy $\sigma$, \Alice{}  chooses 
$ l_n=2$ for all $n\in\w$.
The strategy $\sigma$ instructs \Alice{} to play some cover $\cU_m$.
Then if the set chosen by \Bob{} is of the form $U_{[m,k)}$, then \Alice{} plays the family $\cU_k$.
Each play where \Alice{} uses $\sigma$ has the following form
\[
(2, \cU_{i_0}, U_{[i_0,i_1)}, \cU_{i_1}, U_{[i_1,i_2)}; 2, \cU_{i_2}, U_{[i_2,i_3)}, \cU_{i_3}, U_{[i_3,i_4)};\dotsc)
\]
where $\seq{i_n}{i\in\w}$ is an increasing sequence of natural numbers with $i_0=0$.
Since the sets
\[
a:=\Un_{k\in\w}[i_{2k},i_{2k+1}),\quad b:=\Un_{k\in\w}[i_{2k+1},i_{2k+2})
\]
are disjoint and $a\cup b=\w$, exactly one of them is a member of $X$.
Assume that $a\in X$.
Since $a\notin \bigcup_{k\in\w} U_{[i_{2k+1},i_{2k+2})}$, we have
\[
X\not\sub \Un_{k\in\w}(U_{[i_{2k},i_{2k+1})}\cap U_{[i_{2k+1},i_{2k+2})}).\qedhere
\]
\epf

Combining \cite[Theorem~1.1]{ChoRepZdo15} and \cite[Theorem~10]{Can88}, we 
conclude that $\mathfrak d=\mathfrak c$ implies the existence of a Menger ultrafilter.
\bcor\label{ultraf_d_c}
Assume that $\mathfrak d=\mathfrak c$. There exists a Menger subspace $X$ of 
$2^\w$ such that \Alice{} has a winning strategy in the grouped Menger 
game on $X$, i.e., $X\not\in\mathcal{GM}$.
\ecor

Let us note that there are models of ZFC without Menger ultrafilters.
Indeed, every Menger ultrafilter is a $P$-point, see, e.g., \cite[Observation~3.4]{BelTokZdo16}. According to a result of Shelah published in \cite{Wim82} (see also \cite{ChoGuz}), consistently there are no $P$-points.

\medskip

Next, we introduce a variation of the grouped Menger game which is 
``harder'' for \Alice{} since her choices are more restricted, i.e., 
are not arbitrary open covers of the space in question.
Let $K\sub\Cantor$ be a perfect set.
For $y\in K$, let $\cU_y$ be a countable increasing family of clopen sets in $K$ such that $\Un\cU_y=K\sm\{y\}$.
Let $X\sub K$ be a set such that $K\sm X$ is dense in $K$.
The \emph{weak grouped Menger game played on $X$ in $K$}
($\pBM(K,X)$ in short)
is played as follows: In round $0$ \Alice{} selects a natural number $ l_0>0$, and \Bob{} selects a closed nowhere dense subset $K_0$ of $K$.
Let $L_0:=0$ and $L_1:=l_0$.
Then the players 
play the usual Menger game $ l_0$ subrounds, with the following restrictions. 
In each subround $i\in[L_0,L_1)$,
\Alice{} chooses $y_i\in K\setminus (X\cup K_0)$ and plays the family $\cU_i:=\cU_{y_i}$.
Then \Bob{} replies by choosing 
a set $U_i\in\mathcal U_i$ with $K_0\sub U_i$.

Afterwards, in round $1$ \Alice{} selects a natural number $ l_1>0$, and 
\Bob{} selects a closed nowhere dense set $K_1\sub K$.
Let $L_2:=l_0+l_1$.
Then the players 
play the  Menger game further $ l_1$ subrounds with the restriction given above, i.e., in each subround $i\in[L_1,L_2)$ \Alice{} chooses $y_i\in K\setminus (X\cup K_1)$ and plays the family $\mathcal U_i:=\cU_{y_i}$.
Then \Bob{} replies by choosing 
a set $U_i\in\mathcal U_i$ with $K_1\sub U_i$.

In round $n$ \Alice{} selects a natural number $ l_n>0$, and \Bob{} selects a closed nowhere dense set $K_n\sub K$.
Let $L_{n+1}:=l_0+l_1+\dotsb+l_n$.
Then the players 
play the  Menger game further $ l_n$ subrounds such that in each subround $i\in [L_n,L_{n+1})$ \Alice{} chooses $y_i\in K\setminus (X\cup K_n)$ and plays the family $\mathcal U_i:=\cU_{y_i}$.
Then \Bob{} replies by choosing 
a set $U_i\in\mathcal U_i$ with $K_n\sub U_i$.

\Bob{} wins the game if
\[
X=\bigcup_{n\in\w}\bigcap_{i\in [L_n,L_{n+1})}U_i ,
\]
and \Alice{} wins otherwise.

\brem \label{obs:perf_meag}
Let $X$ be a subset of a perfect set $K\sub\Cantor$ such that $K\sm X$ is dense in $K$. 
If \Alice{} has a winning strategy in the weak grouped Menger game played on $X$ in $K$, then \Alice{} has a winning strategy in the grouped Menger game played on $X$.
\erem

A set $X\sub\Cantor$ is \emph{perfectly meager} if for any perfect set $K\sub \Cantor$, the intersection $X\cap K$ is meager in $K$.

\bprp\label{prp:pmpBM}
Let $K\sub\Cantor$ be a perfect set and $X\sub \Cantor$ be a perfectly meager set.
Then \Bob{} has a winning strategy in the 
weak grouped Menger game played on $X\cap K$ in $K$.
\eprp
\bpf
For each natural number $n$ let $K_n\sub K$ be a closed nowhere dense subset of $K$ such that $X\cap K\sub\Un_{n\in\w}K_n\sub K$. 
Then any strategy for \Bob{} in the weak grouped Menger game played on $X\cap K$ in $K$, where in each round $n$, \Bob{} plays the set $K_n$, is a winning strategy.
\epf

Let $K\sub \Cantor$ be a perfect set and $X\sub\Cantor$.
The diagram below presents the relations between the considered properties.
By $\mathrm{A}\not\uparrow \pBM(K,X\cap K)$ we mean that \Alice{} has no winning strategy and by $\mathrm{B}\uparrow \pBM(K,X\cap K)$ we mean that \Bob{} has a winning strategy in the game $\pBM(K,X\cap K)$.
Note that co-consonant spaces are preserved by closed subspaces.

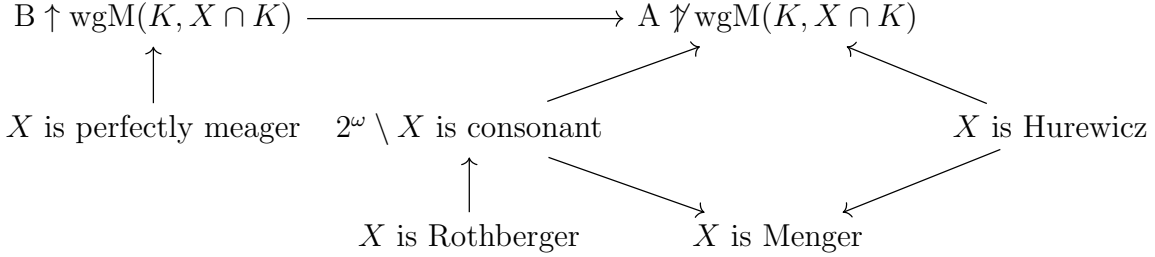
\begin{figure}[H]
{
\begin{tikzcd}[ampersand replacement=\&,column sep=.1cm]
\textrm{B}\uparrow \pBM(K,X\cap K)\arrow{rr}\& {}\&\textrm{A}\not\uparrow \pBM(K,X\cap K)\&{}\\
X\text{ is perfectly meager}\arrow{u}\& \Cantor\sm X\text{ is consonant}\arrow{rd}\arrow{ru}\& {} \& X\text{ is Hurewicz}\arrow{ld}\arrow{lu}\\
{}\&{X\text{ is Rothberger}}\arrow{u}\& X\text{ is Menger}\&{}
\end{tikzcd}
}    
\end{figure}

The next fact is similar to Lemma~\ref{clm:easy}.

\blem \label{clm:easy_long}
We use the notation and objects described in Section~2.
Suppose that
$X$ is a subset of the perfect set $K\sub(\Cantor)^S$ such that  
$K\sm X$ is dense in $K$ and \Alice{}
has no winning strategy in $\pBM(K,X)$.

Then there exists a sequence 
$\seq{\la i_n,j_n,C_n\ra}{n\in\w}$  such that
\begin{enumerate}
\item $\seq{i_n}{n\in\w}$ is a strictly increasing sequence of natural numbers;
\item $C_n\sub\Sigma_{i_n}$; 
\item for every $n\in\w$, we have $j_n\in [i_n,i_{n+1})$,
and for every $\nu\in C_n$ there exists $e_n(\nu)\in \Sigma_{j_n}$ extending
$\nu$, such that 
$$C_{n+1}=\bigcup_{\nu\in C_n}\{\sigma\in\Sigma_{i_{n+1}}\: :\: \sigma \succ e_n(\nu)\};$$
\item the maps $e_n\colon C_n\to \Sigma_{j_n}$ are coherent; and
\item $\bigcap_{n\in\w}\bigcup_{\nu\in C_n}[\nu]\cap X=\emptyset$.
\end{enumerate}
\elem

\bpf
We shall describe a strategy $\S$ for \Alice{} in the weak grouped Menger game played on $X$ in $K$ 
such that
each play lost \Alice{} gives rise to the objects whose existence we need to establish.

\emph{Round 0.} Let $C_0:=\Sigma_0$ and $i_0:=0$.
\Alice{} declares that the $0$-th group will have 
length $ l_0:=|C_0|$. Let $\{\nu_j:j< l_0\}\sub\Sigma_0$ be an enumeration of
$C_0$.
Suppose that
\Bob{} plays a closed nowhere dense $K_0\sub K$.

\emph{Subround $\la 0,0\ra$.}
By density of $K\sm X$, \Alice{} picks $y_{\la 0,0\ra}\in [\nu_0]\setminus (X\cup K_0)$ and she plays the family $\mathcal U_{\la 0,0\ra}:=\cU_{y_{\la 0,0\ra}}$.
Suppose that 
\Bob{} replies by choosing $U_{\la 0,0\ra}\in\mathcal U_{\la 0,0\ra}$
with $K_0\sub  U_{\la 0,0\ra}$.
Take $j_{\la 0,0\ra}>i_0$ and $\nu_{\la 0,0\ra}\in\Sigma_{j_{\la 0,0\ra}}$ with $\nu_{\la 0,0\ra}\succ\nu_0$ such that $y_{\la 0,0\ra}\in [\nu_{\la 0,0\ra}]$
and $[\nu_{\la 0,0\ra}]\cap U_{\la 0,0\ra}=\emptyset$.
Let $e_{\la 0,0\ra}\colon C_0\to\Sigma_{j_{\la 0,0\ra}}$
be a coherent map such that $e_{\la 0,0\ra}(\nu_0)=\nu_{\la 0,0\ra}$.
Fix a natural number $a$ with $0<a< l_0$ and assume that the players have already defined the following sequences:
\begin{itemize}
\item $\seq{y_{\la 0,b\ra}}{0\leq b<a}$ of elements $y_{\la 0,b\ra}\in K\setminus (X\cup K_0)$,
\item $\seq{\mathcal U_{\la 0,b\ra}}{0\leq b<a}$ of covers of $X$ by clopen subsets of $K$ such that $\cU_{\la 0, b\ra}=\cU_{y_{\la 0, b\ra}}$,
\item
$\seq{U_{\la 0,b\ra}}{0\leq b<a}$ of sets $U_{\la 0,b\ra}\in
\mathcal U_{\la 0,b\ra}$ with $K_0\sub U_{\la 0,b\ra}$,
\item
$\seq{j_{\la 0,b\ra}}{0\leq b<a}\in\w^{\uparrow a}$ with $i_0<j_{\la 0,0\ra}$,
\item $\seq{\nu_{\la 0,b\ra}}{0\leq b<a}$ of maps $\nu_{\la 0,b\ra}\in\Sigma_{j_{\la 0,b\ra}}$ such that $y_{\la 0,b\ra}\in [\nu_{\la 0,b\ra}]$, $\nu_b\prec\nu_{\la 0,b\ra}$ and $[\nu_{\la 0,b\ra}]\cap U_{\la 0,b\ra}=\emptyset$,
\item
 $\seq{e_{\la 0,b\ra}}{0\leq b<a}$ of coherent maps $e_{\la 0,b\ra}\colon C_0\to\Sigma_{j_{\la 0,b\ra}}$ with  $e_{\la 0,b\ra}(\nu_b)=\nu_{\la 0,b\ra}$ and for every
$0\leq b_0<b_1<a$  and $j< l_0$ we have $e_{\la 0,b_0\ra}(\nu_{j})\prec e_{\la 0,b_1\ra}(\nu_{j})$.
\end{itemize}

\emph{Subround $\la 0,a\ra$}:
\Alice{} picks $y_{\la 0,a\ra}\in [e_{\la 0,a-1\ra}(\nu_a)]\setminus (X\cup K_0)$ and she plays the family $\mathcal U_{\la 0,a\ra}:=\cU_{y_{\la 0, a\ra}}$.
Suppose that 
\Bob{} replies by choosing $U_{\la 0,a\ra}\in\mathcal U_{\la 0,a\ra}$
with $K_0\sub  U_{\la 0,a\ra}$.
Take $j_{\la 0,a\ra}>j_{\la 0,a-1\ra}$ and $\nu_{\la 0,a\ra}\in\Sigma_{j_{\la 0,a\ra}}$ with $\nu_{\la 0,a\ra}\succ e_{\la 0,a-1\ra}(\nu_a)$ such that $y_{\la 0,a\ra}\in [\nu_{\la 0,a\ra}]$
and $[\nu_{\la 0,a\ra}]\cap U_{\la 0,a\ra}=\emptyset$.
Let $e'_{\la 0,a\ra}\colon e_{\la 0,a-1\ra}[ C_0]\to\Sigma_{j_{\la 0,a\ra}}$
be a coherent map such that $e'_{\la 0,a\ra}(e_{\la 0,a-1\ra}(\nu_a))=\nu_{\la 0,a\ra}$.
Put $e_{\la 0,a\ra}:=e'_{\la 0,a\ra}\circ e_{\la 0,a-1\ra}$.

After subround $\la 0, l_0-1\ra$, the last subround of round 0, we set 
\[
j_0:=j_{\la 0, l_0-1\ra},\quad i_1:=j_0+1, \quad e_0:=e_{\la 0, l_0-1\ra},\quad C_1:=\bigcup_{\nu\in C_0}\sset{\sigma\in\Sigma_{i_1}}{e_0(\nu)\prec\sigma}.
\]

Fix a natural number $n>0$ and assume that the elements of the sequence
\[
\la i_0,C_0,j_0,e_0;\dotsc;i_{n-1},C_{n-1},j_{n-1},e_{n-1}; i_{n},C_{n}\ra
\]
satisfy all relevant instances of $(1)$--$(4)$.

\emph{Round $n$}.
\Alice{} declares that the $n$-th group will have 
length $ l_{n}:=|C_{n}|$.
Let $\{\nu_j:j< l_{n}\}\sub\Sigma_{i_{n}}$ be an enumeration\footnote{Formally, we should have written $\sset{\nu^{n}_j}{j< l_{n}}$ instead of $\sset{\nu_j}{j< l_{n}}$, but
we   omit extra indices in order to shorten our notation.}
of
$C_{n}$. Suppose that \Bob{} plays a closed nowhere dense subset $K_{n}\sub K$.

\emph{Subround $\la n,0\ra$}:
\Alice{} picks $y_{\la n,0\ra}\in [\nu_0]\setminus (X\cup K_{n})$ and she plays the family $\mathcal U_{\la n,0\ra}:=\cU_{y_{\la n,0\ra}}$.
Suppose that 
\Bob{} replies by choosing $U_{\la n,0\ra}\in\mathcal U_{\la n,0\ra}$ with  $K_n\sub  U_{\la n, 0\ra}$.
Take $j_{\la n,0\ra}>i_{n}$ and $\nu_{\la n,0\ra}\in\Sigma_{j_{\la n,0\ra}}$ with $\nu_{\la n,0\ra}\succ\nu_0$ such that $y_{\la n,0\ra}\in [\nu_{\la n,0\ra}]$
and $[\nu_{\la n,0\ra}]\cap U_{\la n,0\ra}=\emptyset$.
Let $e_{\la n,0\ra}\colon C_{n}\to\Sigma_{j_{\la n,0\ra}}$
be a coherent map such that $e_{\la n,0\ra}(\nu_0)=\nu_{\la n,0\ra}$.

Fix a natural number $a$ with $0<a< l_{n}$  and assume that the players have already defined the following sequences:
\begin{itemize}
\item $\seq{y_{\la n,b\ra}}{0\leq b<a}$ of elements $y_{\la n,b\ra}\in K\setminus (X\cup K_{n})$;
\item 
$\seq{\mathcal U_{\la n,b\ra}}{0\leq b<a}$ of covers of $X$ by clopen subsets of $K$ such that $\mathcal U_{\la n,b\ra}=\cU_{y_{\la n,b\ra}}$;
\item
$\seq{U_{\la n,b\ra}}{0\leq b<a}$ of sets $U_{\la n,b\ra}\in
\mathcal U_{\la n,b\ra}$ with $K_{n}\sub U_{\la n,b\ra}$;
\item
$\seq{j_{\la n,b\ra}}{0\leq b<a}$ which is an increasing sequence with $j_{\la n,0\ra}>i_{n}$;
\item $\seq{\nu_{\la n,b\ra}}{0\leq b<a}$ of maps $\nu_{\la n,b\ra}\in\Sigma_{j_{\la n,b\ra}}$ such that $y_{\la n,b\ra}\in [\nu_{\la n,b\ra}]$, $\nu_b\prec\nu_{\la n,b\ra}$ and $[\nu_{\la n,b\ra}]\cap U_{\la n,b\ra}=\emptyset$;
\item
 $\seq{e_{\la n,b\ra}}{0\leq b<a}$ of coherent maps $e_{\la n,b\ra}\colon C_{n}\to\Sigma_{j_{\la n,b\ra}}$ with $e_{\la n,b\ra}(\nu_b)=\nu_{\la n,b\ra}$ and for every
$0\leq b_0<b_1<a$  and $j< l_{n}$ we have $e_{\la n,b_0\ra}(\nu_{j})\prec e_{\la n,b_1\ra}(\nu_{j})$.
\end{itemize}

\emph{Subround $\la n,a\ra$:}
\Alice{} picks $y_{\la n,a\ra}\in [e_{\la n,a-1\ra}(\nu_a)]\setminus (X\cup K_{n})$ and she plays the family $\mathcal U_{\la n,a\ra}:=\cU_{y_{\la n,a\ra}}$.
Suppose that 
\Bob{} replies by choosing $U_{\la n,a\ra}\in\mathcal U_{\la n,a\ra}$ such that
$K_n\sub  U_{\la n,a\ra}$.
Take $j_{\la n,a\ra}>j_{\la n,a-1\ra}$ and $\nu_{\la n,a\ra}\in\Sigma_{j_{\la n,a\ra}}$ with $\nu_{\la n,a\ra}\succ e_{\la n,a-1\ra}(\nu_a)$ such that $y_{\la n,a\ra}\in [\nu_{\la n,a\ra}]$
and $[\nu_{\la n,a\ra}]\cap U_{\la n,a\ra}=\emptyset$.
Let $e'_{\la n,a\ra}\colon e_{\la n,a-1\ra}[ C_{n}]\to\Sigma_{j_{\la n,a\ra}}$
be a coherent map such that $e'_{\la n,a\ra}(e_{\la n,a-1\ra}(\nu_a))=\nu_{\la n,a\ra}$.
Put $e_{\la n,a\ra}:=e'_{\la n,a\ra}\circ e_{\la n,a-1\ra}$.

After subround $\la n, l_{n}-1\ra$, the last subround of round $n$, we set 
\[
\begin{gathered}
j_{n}:=j_{\la n, l_{n}-1\ra},\quad
i_{n+1}:=j_{n}+1,\quad
e_{n}:=e_{\la n, l_{n}-1\ra},\\
C_{n+1}:=\bigcup_{\nu\in C_{n}}\{\sigma\in\Sigma_{i_{n+1}}:e_{n}(\nu)\prec\sigma\}.
\end{gathered}
\]

This completes the definition of the strategy $\S$ in the weak grouped Menger game played on $X$. 
It remains to notice that 
any play in this game in which \Alice{} uses $\S$ gives a sequence of objects we require in our lemma, such that conditions $(1)$--$(4)$ are satisfied, and
if \Alice{} loses (and such a play exists because $\S$ cannot be winning), then
also $(5)$ is satisfied.
\epf


\section{The weak grouped Menger game and Sacks forcing } \label{menger_size_rev}

Again, by $V$ we mean a ground model of ZFC and $G_{\w_2}$ is an $\mathbb S_{\w_2}$-generic filter over $V$.

\bthm\label{thm:main_new}\mbox{}
Assume that $V$ satisfies CH.
In $V[G_{\w_2}]$,  suppose that $X\sub 2^\w$ and for any perfect set $K\sub 2^\w$ such that $K\setminus X$ is dense in $K$ \Alice{} has no winning strategy in the weak grouped Menger game played on $X\cap K$ in $K$. 
Then both $X$ and $2^\w\setminus X$ are unions of $\w_1$ compact sets.
 \ethm

For the proof of Theorem~\ref{thm:main_new}, we need the following auxiliary result.
Similarly to Lemma~\ref{lem:LS_ti}, it can be proved in the same way as \cite[Lemma 5.10]{BlaShe87},
a rather standard  argument is left to the reader.
\blem\label{lem:LS}
In $V[G_{\w_2}]$, let $X\sub\Cantor$. Then there exists a limit ordinal $\alpha<\w_2$ of cofinality $\w_1$
such that 
\be
\item $X\cap V[G_\alpha]\in V[G_\alpha]$, and if $K,K'\sub 2^\w$ are closed 
crowded sets and coded in $V[G_\alpha]$, and  $K'\sub K\setminus (X\cap V[G_\alpha])$,
then $K'\sub K\setminus X$;
\item There is a function in $V[G_\alpha]$ which assigns to every perfect set $K\sub\Cantor$, coded in $V[G_\alpha]$, such that $K\sm X$ is not dense in $K$, a nonempty clopen subset $O$ of $2^\w$ with\footnote{This function assigns to each code in $V[G_\alpha]$ for a perfect set a code in $V[G_\alpha]$ for a basic open set with the given properties.} $\emptyset\neq K\cap O\sub K\cap X$.
\ee
Moreover, if for each perfect set $K\sub \Cantor$ such that $K\sm X$ is dense in $K$,  \Alice{} has no winning strategy in the weak grouped Menger game played on $X\cap K$ in $K$,
then we can in addition assume that
\begin{itemize} 
\item[(3)]  in $V[G_\alpha]$,  
for each perfect set $K\sub \Cantor$ such that $K\sm (X\cap V[G_\alpha])$ is dense in $K$,  \Alice{} has no winning strategy in the weak grouped Menger game played on $X\cap K\cap V[G_\alpha]$ in $K$.
\end{itemize}
\elem


\bpf[{Proof of Theorem~\ref{thm:main_new}}]
Let $\alpha$ be such as in Lemma~\ref{lem:LS}. Working in $V[G_{\w_2}]$,
we claim that 
\begin{equation} \label{eq1}
X=\bigcup\big\{L\sub 2^\w\: :\: L\sub X,\mbox{ $L$ is compact, and $L$ is coded in }V[G_\alpha]\big\} \mbox{ and}
\end{equation}
\begin{equation} \label{eq2}
2^\w\setminus X=\bigcup\big\{L\sub 2^\w\: :\: L\sub 2^\w\setminus X,\mbox{ $L$ is compact, and $L$ is coded in }V[G_\alpha]\big\}.
\end{equation}
Since in $V[G_\alpha]$ the remainder 
$\mathbb S_{\alpha,\w_2}$ is order-isomorphic to $\mathbb S_{\w_2}$, there is no loss of generality in  
assuming that $\alpha=0$, i.e., that $V=V[G_\alpha]$. 

First we prove (\ref{eq1}). 
Let us pick $z\in X$ and let $\gamma$ be the minimal ordinal with 
$z\in V[G_\gamma]$.
From now on, whenever we do not specify that we work in some other model, we work in $V$.
Let $r\in G_{\w_2}$ and $\tau\in V^{\mathbb S_\gamma}$ be 
such that $\tau^{G_\gamma}=z$ and $r\Vdash\tau\in X$. 
Let $p_0 \geq r \restriction \gamma$ be 
such that $p_0\in G_{\gamma}$ and
$$p_0\Vdash_\gamma\tau \in (2^\w \cap V[\dot{G}_\gamma])\setminus\bigcup_{\beta<\gamma}V[\dot{G}_\beta].$$ 
We shall find $q\geq p_0$, $q \in \mathbb{S}_\gamma$ and a compact set $L\sub X$
coded in $V$ such that $q\Vdash_{\w_2}\tau\in L$. This would accomplish the proof: The genericity
of $G_{\gamma}$ implies that there is a condition $q$ as above which lies in $G_\gamma \subseteq G_{\w_2}$.
Then $z = \tau^{G_\gamma} = \tau^{G_{\w_2}}\in L$. 

Take $p$ and $F_n,k_n,l_n, y_{\sigma_n}$ from Lemma~\ref{lem:Mill}, applied to $p_0$ and  $\tau$.
Let $\Sigma_n$ be the set of all maps $\sigma\colon F_n\to 2^{k_n}$ consistent with $p$, where  $n\in\w$.
By Lemma~\ref{lem:prompt}(1), for each map $\sigma\in\Sigma_n$, there are maps $\sigma'\neq \sigma''\in\Sigma_{n+1}$ extending $\sigma$, which implies $[\sigma']\cap[\sigma'']=\emptyset$.
For $S:=\supp(p)$, define a perfect set $K$, exactly in the same way as before Lemma~\ref{coh_sel}, i.e., 
\[
K:=\bigcap_{n\in\w}\Un\sset{[\sigma]}{\sigma\in\Sigma_n}.
\]
Then the family of all sets $[\sigma]$, where $\sigma\in\bigcup_{n\in\w}\Sigma_n$, is a basis for $K$.
Note also that Lemma~\ref{lem:prompt}(3) ensures that 
$S$ and the sequence $\seq{k_n,F_n,\Sigma_n}{n\in\w}$
satisfy  $(e_f)$ stated at the beginning of Section~\ref{sec2}.

Let $h:K\to 2^\w$ be defined in the same way as in the proof of Theorem~\ref{thm:main_ti}.
Thus $h\colon K\to h[K]$ is a homeomorphism, and hence $h[K]$ is perfect.
Fix a natural number $n$.
By Lemma~\ref{lem:Mill}(2) and Observation~\ref{obs:det2}(3), the set $\sset{p|\sigma_n}{\sigma_n\in\Sigma_n}$  is a maximal antichain above $p$.
Applying Lemma~\ref{lem:Mill}(4), we have that 
\[
p\forces \tau\restrict l_n\in \sset{y_{\sigma_n}}{\sigma_n\in\Sigma_n}.
\]
It follows from the above and from the definition of the function $h$  that
$p\forces \tau\in h[K].$

Assume that in $V[G_{\w_2}]$, the set $h[K]\setminus X$ is dense in $h[K]$.
We shall show that this is impossible.
By the assumption, \Alice{} has no winning strategy in $\pBM(h[K], X\cap h[K])$.
We proceed in $V$. 
By Lemma~\ref{lem:LS}(3) we have that 
\Alice{} has no winning strategy in $\pBM(h[K],h[K]\cap X\cap V)$.
Since $h$ is a homeomorphism and $h\inv[X\cap V]=h\inv[X]\cap V=h\inv[X]\cap V\cap K$ (because $h$ is defined in $V$ and its domain is $K$), \Alice{} has no winning strategy in $\pBM(K,h\inv[X]\cap V)$. Applying Lemma~\ref{clm:easy_long} to
$S=\mathrm{supp}(p)$ and $h^{-1}[X]\cap V\sub K$, we can get 
a sequence $\seq{\la i_n,j_n,C_n\ra}{n\in\w}$ satisfying the conclusion of
Lemma~\ref{clm:easy_long}. 
Lemma~\ref{clm:easy_long}(5) yields
\[
K':=\bigcap_{n\in\w}\Un\sset{[\nu]}{\nu\in C_n}\sub K\sm (h\inv[X]\cap V),
\]
and therefore $h[K']\sub  h[K]\setminus (X\cap V)$, both of these inclusions holding in $V$.
Applying Lemma~\ref{lem:LS}(1) we conclude that 
$h[K']\sub  h[K]\setminus X$  holds in $V[G_{\w_2}]$.

 Let $q\geq p$ be a condition given by Lemma~\ref{clm:Sigma}.
 Note that $q\Vdash\tau \restriction l_{i_n} \in \sset{y_{\nu}}{\nu\in C_n}$ because
$\{p|\nu:\nu\in C_n\}$ is predense above $q$, for all $n\in\w$, which implies
$q\Vdash\tau\in h[K']$. It follows from the above that
$q\Vdash \tau\not\in X$, which is impossible since $q^\frown r \restriction [\gamma,\w_2) \Vdash \tau \in X$. 

Now assume that in $V[G_{\w_2}]$, we have $\Int_{h[K]}(h[K]\cap X)\neq\emptyset$.
Since $h[K]$ is coded in $V$, there is a clopen set $O$  (namely the one assigned to $h[K]$
by the function given by Lemma~\ref{lem:LS}(2))
such that  $\emptyset\neq h[K]\cap O\sub X$. 
It follows from the definition of $h$ that there exist $n\in\w$ and $\sigma\in \Sigma_n$
with $h[[\sigma]\cap K]\sub O$, and hence $h[[\sigma]\cap K]\sub X$. 
Let $q:=p|\sigma$.
Since $q\Vdash\tau\in h[[\sigma]\cap K]$, the set $L:=h[[\sigma]\cap K]$
is a perfect set coded in $V$ which is a subset of $X$
and $q\Vdash \tau\in L$. This completes the proof of Equation~\ref{eq1}.

Next, we prove (\ref{eq2}). The argument is very similar to that of
 (\ref{eq1}), but we anyway give it for the sake of completeness. 
As above, let $r\in G_{\w_2}$, $p_0 \geq r \restriction \gamma$, $p_0\in G_{\gamma}$ and $\tau\in V^{\mathbb S_\gamma}$ be 
such that $r\Vdash\tau\in 2^\w \setminus X$ and 
$$p_0\Vdash_\gamma\tau \in (2^\w \cap V[\dot{G}_\gamma])\setminus\bigcup_{\beta<\gamma}V[\dot{G}_\beta].$$
We shall find $q\geq p_0$, $q \in \mathbb{S}_\gamma$ and a compact set $L\sub 2^\w\setminus X$
coded in $V$ such that $q\Vdash_{\w_2} \tau\in L$.
As in the case of (\ref{eq1}), this would accomplish the proof. 

Again, take $p$ and $F_n,k_n,l_n, y_{\sigma_n}$ from Lemma~\ref{lem:Mill}, applied to $p_0$ and  $\tau$.
Let $\Sigma_n$, $[\sigma]$ for $\sigma\in\Sigma_n$,
$K$, and $h:K\to 2^\w$ be defined in the same way as in the proof of (\ref{eq1}).
We have
$p\forces \tau\in h[K].$ 

Assume that in $V[G_{\w_2}]$, we have $\Int_{h[K]}(h[K]\cap X)\neq\emptyset$.
We shall show that this case is impossible.
Since $h[K]$ is coded in $V$, there is a clopen set $O$  
such that  $\emptyset\neq h[K]\cap O\sub X$. 
It follows from the above that there exist $n\in\w$ and $\sigma\in \Sigma_n$
with $h[[\sigma]\cap K]\sub O$, and hence $h[[\sigma]\cap K]\sub X$. 
Let $q:=p|\sigma$.
Since $q\Vdash\tau\in h[[\sigma]\cap K]$, the set $L:=h[[\sigma]\cap K]$
is a perfect set coded in $V$ which is a subset of $X$ 
and $q\Vdash \tau\in L$. Consequently, $q\Vdash \tau\in X$, a contradiction to
$r\Vdash\tau\in 2^\w\setminus X$.

Now assume that in  $V[G_{\w_2}]$, the set $h[K]\setminus X$ is dense in $h[K]$.
In $V$, we have that
$h[K]\setminus (X\cap V)$ is dense in $h[K]$: If there were a clopen $K'\sub  2^\w$
with 
\[
\emptyset\neq K'\cap h[K]\sub  h[K]\setminus (X\cap V),
\]
then we would get $K'\cap h[K]\sub  h[K]\setminus X$ holding in $V[G_{\w_2}]$
by Lemma~\ref{lem:LS}(1).

By the assumption, \Alice{} has no winning strategy in the game $\pBM(h[K],X\cap h[K])$. 
By Lemma~\ref{lem:LS}(3), in $V,$ \Alice{} has no winning strategy in  $\pBM(h[K],X\cap h[K]\cap V)$.
Since $h$ is a homeomorphism, \Alice{} has no winning strategy in $\pBM(K,h\inv[X]\cap K\cap V)$.
Now let $\seq{\la i_n,j_n,C_n\ra}{n\in\w}$, 
$K'$ and $q$ be the same as in the proof of (\ref{eq1}). Again, since
$K'\sub  K\setminus(h\inv[X]\cap V)$ holds in $V$, we have
$K'\sub  K\setminus h\inv[X]$ in $V[G_{\w_2}]$ by  Lemma~\ref{lem:LS}(1), or equivalently $h[K']\sub  h[K]\setminus X$.
Repeating our previous arguments we get
\[q\Vdash\tau\in h[K']\sub 2^\w\setminus X.\] It follows from the above that the set 
 $L:=h[K']$ is a perfect subset of $2^\w\setminus X$ and
$q\Vdash \tau\in L$, which  completes the proof of (\ref{eq2}).
\epf

Combining Theorem~\ref{thm:main_new} with Proposition~\ref{obs_obs}, Remark~\ref{obs:perf_meag} and Proposition~\ref{prp:pmpBM},
we get the following result.

\bcor\label{cor_resqued}
Assume that $V$ satisfies GCH.
 In $V[G_{\w_2}]$, the following assertions hold.
\be 
\item Each consonant (Hurewicz) subset of $\Cantor$ as well as its complement are unions of $\fd=\w_1$ many compact subspaces. In particular, there are
$\mathfrak c=\w_2$ consonant (Hurewicz) subspaces of $2^\w$.
\item Each perfectly meager subset of $\Cantor$ has cardinality at most $\fd=\w_1$, and its complement is a union of $\w_1$ compact sets.
\ee
\ecor
The first half of Corollary~\ref{cor_resqued}(2), namely that
all perfectly meager subsets of $\Cantor$ have cardinality at most $\w_1$
in the Sacks model, was established by Miller~\cite[\S~5]{mill}.

\brem\label{rem}
In $V[G_{\w_2}]$, considered in the above corollary, there is a Luzin subset of $\Cantor$, i.e., an uncountable set whose intersection with any meager set is at most countable, which is totally imperfect and Menger~\cite{Hu27}, but it is not (perfectly) meager.
There is also a perfectly meager set that is not Menger:
Since $\fd=\w_1$ in $V[G_{\w_2}]$, there is a dominating set $X=\sset{x_\alpha}{\alpha<\w_1}$ in $\roth$, where $x_\beta\les x_\alpha$ for all ordinal numbers $\beta<\alpha<\w_1$.
This set is not Menger.
The set $X\cup\Fin$ satisfies the Hurewicz covering property \cite[Theorem~10.]{BaTs}, and any totally imperfect set with this property is perfectly meager~\cite[Theorem~5.5.]{coc2}.
Thus, the set $X$ is perfectly meager, too. Alternatively, we could use here the main result of \cite{Plewik}, which implies directly that $X$ is perfectly meager since $\leq^*$
is a Borel subset of  $\w^\w\times \w^\w$.

\erem

Theorems~\ref{thm:main_ti} and \ref{thm:main_new} motivate the following problem.

\bprb
\label{Menger_full}
In the Sacks model:
\begin{enumerate}

\item  Is every Menger space $X\sub 2^\w$ a union of $\w_1$-many of its compact subspaces?
\item Is the complement $2^\w\setminus X$ of a  Menger set $X\sub 2^\w$ a union of $\w_1$-many of its compact subspaces?
\item  Are there only $\w_2$-many Menger subsets of $2^\w$?
\item  Is the complement $2^\w\setminus X$ of a totally imperfect Menger set $X\sub 2^\w$ a union of $\w_1$-many of its compact subspaces?
\end{enumerate}
\eprb

Regarding the last item of Problem~\ref{Menger_full},
we do not know the answer even to the following question.
\bprb \label{8_6}
In the Sacks model,
is the complement $2^\w\setminus X$ of any set $X\sub\Cantor$ with $\card{X}=\w_1$ a union of $\w_1$-many of its compact subspaces? In particular, 
is $2^\w\setminus(2^\w\cap V)$ a union of
$\w_1$-many of its compact subspaces?
\eprb

\section{Menger sets and Hechler forcing}\label{Hech}

The results from the previous section lead to the question, whether $\fd<\fc$ implies that any totally imperfect Menger subset of $\Cantor$ has cardinality at most $\fd$.
We address this problem, showing that this is not the case.
We also provide a consistent result that the size of the family of all Hurewicz subsets of $\Cantor$ can be equal $2^\fc$ even if $\fd<\fc$.
Let $\bbP$ be a definable ccc forcing notion of size $\mathfrak c$ which adds dominating reals over a ground model (e.g., the poset defined in \cite[p.~95]{BauDor85}, nowadays commonly named  Hechler forcing)  and $\bbP_{\w_1}$ be an iterated forcing of length $\w_1$ with finite support, where each iterand is equal to $\bbP$.
In this section, by $G_{\w_1}$ we mean a $\bbP_{\w_1}$-generic filter over a ground model $V$ of ZFC.

\bprp\label{prp:hech}
In $V[G_{\w_1}]$, each subset of $\Cantor\cap V$ is Hurewicz.
\eprp

\bpf
Let $X\sub \Cantor\cap V$ and $\phi\colon X\to\roth$ be a continuous function.
The function $\phi$ can be extended to a continuous function $\Phi$ defined on a $G_\delta$-set $A$, containing $X$.
Then there is an ordinal number $\alpha<\w_1$ such that the function $\Phi$ and the set $A$ are coded in $V[G_\alpha]$.
It follows that $\Phi[X]\sub V[G_\alpha]$. 
Since there exists  a function $g_\alpha$, added in step $\alpha$, which dominates any real from $V[G_\alpha]$, the set $\Phi[X]$ is bounded in $V[G_{\w_1}]$.
By the Hurewicz--Rec{\l}aw characterization of the Hurewicz property~\cite[Proposition~1]{reclaw}, the set $X$ is Hurewicz\footnote{As noted by the referee, a similar argument actually gives that
$\phi[X]$ is bounded for any Borel function $f:X\to\w^\w$.} in $V[G_{\w_1}]$.
\epf

\bthm
Assume that $V$ satisfies $\neg$CH.
In $V[G_{\w_1}]$, we have $\fd<\fc$ and the following assertions hold.
\be 
\item
There is a totally imperfect Hurewicz and Rothberger (and thus  Menger) subset of $\Cantor$ with cardinality $\fc$.
\item 
There are $2^\fc$-many Hurewicz and Rothberger (and thus consonant) subsets of $\Cantor$.
\ee
\ethm

\bpf
Since $\bbP_{\w_1}$ adds $\w_1$ dominating reals to the ground model, we have $\fd=\w_1<\fc$ in $V[G_{\w_1}]$.
Since the forcing $\bbP$ is ccc, the forcing $\bbP_{\w_1}$ is ccc, too. 
Since $\bbP_{\w_1}$ adds reals to the ground model, any subset of $\Cantor\cap V$ is totally imperfect in $V[G_{\w_1}]$.
By Proposition~\ref{prp:hech}, any subset of $\Cantor\cap V$ is Hurewicz in $V[G_{\w_1}]$.

Since Cohen reals are added by any tail of the considered iteration, any subset of 
$2^\w\cap V$ is also Rothberger in $V[G_{\w_1}]$, by virtue of an argument  similar to that in the proof 
of Proposition~\ref{prp:hech}, more details could be found in the proof of \cite[Theorem~11]{SchTal10}.
\epf

\section{Comments and open problems}

Let $\PN$ be the power set of $\w$, the set of natural numbers.
We identify each element of $\PN$ with its characteristic function, an element of $\Cantor$.
In that way we introduce a topology on $\PN$.
Let $\roth$ be the family of all infinite sets in $\PN$.
Each set in $\roth$ we identify with an increasing enumeration of its elements, a function in the Baire space $\NN$.
We have $\roth\sub\NN$ and topologies in $\roth$ induced from $\PN$ and $\NN$ are the same.
Let $\Fin$ be the family of all finite sets in $\PN$.
A totally imperfect Menger set constructed by Bartoszy\'{n}ski and Tsaban, mentioned in Theorem~\ref{thm:bats} is a set the form $X\cup\Fin\sub\PN$ of size $\fd$, such that for any function $d\in\roth$, we have $\card{\sset{x\in X}{x\les d }}<\fd$.
In fact, any set with these properties is totally imperfect and Menger~\cite[Remark~18.]{BaTs}.
It has also a stronger covering property $\sgo$, described in details in Subsection~8.2.

\bprp
There are at least $2^\fd$-many totally imperfect Menger subsets of $\Cantor$.
\eprp
\bpf
Let $X\sub \roth$ be a set of size $\fd$ such that all coordinates of elements in $X$ are even numbers and for any function $d\in\roth$, we have $\card{\sset{x\in X}{x\les d }}<\fd$.
Let $\sset{x_\alpha}{\alpha<\fd}$ be a bijective enumeration of elements in $X$.
Fix a function $y\colon\fd\to\Cantor$.
For each ordinal number $\alpha<\fd$, let $x_\alpha+y_\alpha$ be the function such that $(x_\alpha+y_\alpha)(n):=x_\alpha(n)+y_\alpha(n)$ for all $n$.
Then $x_\alpha+y_\alpha\in\roth$ and $x_\alpha\les x_\alpha+y_\alpha$ for all ordinal numbers $\alpha<\fd$.
Thus, for any function $d\in\roth$, we have
$\sset{\alpha}{x_\alpha+y_\alpha\les d}\sub \sset{\alpha}{x_\alpha\les d}$, and the latter set has size smaller than $\fd$.
Then the set $\sset{x_\alpha+y_\alpha}{\alpha<\fd}\cup\Fin$ is totally imperfect and Menger.
It follows that for different functions $y\colon\fd\to\Cantor$, we get different totally imperfect Menger sets.
\epf

\subsection{Menger vs. property $\sgo$}

A cover of a space is a $\gamma$-cover if it is infinite and any point of the space belongs to all but finitely many sets in the cover. 
A space satisfies property $\sgo$ if for any sequence $\eseq{\cU}$ of open $\gamma$-covers of the space there are sets $U_0\in\cU_0, U_1\in\cU_1,\dotsc$ such that the family $\sset{U_n}{n\in\w}$ covers the space.
This property implies the Menger property.
By the result of Just, Miller, Scheepers and Szeptycki~\cite[Theorem~2.3.]{coc2} (see also the work of Sakai~\cite[Lemma~2.1.]{sakai}), any subset of $\Cantor$ satisfying $\sgo$ is totally imperfect.
The following questions are one of the major open problems in the combinatorial covering properties theory.

\bprb\mbox{} \label{8_7}
\be
\item
Is there a ZFC example of a totally imperfect Menger subset of $\Cantor$ which does not satisfy $\sgo$?
\item Is there a subset of $\Cantor$ whose continuous images into $\Cantor$ are totally imperfect and Menger, which does not satisfy $\sgo$?
\ee
\eprb
\noindent In the first item of the problem above we ask about ZFC examples because
 under CH there exists even a Hurewicz totally imperfect 
subspace of $2^\omega$ which can be mapped continuously onto $2^\omega$, see \cite{SzeWaiZdo??}. 
The second item of Problem~\ref{8_7}  is motivated by the fact that the property $\sgo$ is preserved by continuous functions.
By the results from Section~2, we can put this problem in a more specific context.

\bprb
Let $G_{\w_2}$ ba an $\Pwtwo$-generic filter over a ground model $V$.
In $V[G_{\w_2}]$, does any Menger set of cardinality $\w_1$ satisfy $\sgo$?
\eprb

\subsection{Other problems}

The following problem is motivated by Remark~\ref{rem}.
\bprb
Is any perfectly meager subset of $\Cantor$ contained in a Menger totally imperfect set?
\eprb

We do not know whether the conclusion of 
Corollary~\ref{ultraf_d_c} holds in ZFC.
\bprb\label{alice_yfc}
Is it consistent that $\mathcal{GM}$ coincides with the family of all 
Menger subspaces of $2^\w$? In other words, is it consistent that for every Menger $X\sub 2^\w$ \Alice{} has no winning strategy in the grouped Menger game on $X$?
\eprb

\end{document}

%% file: shortcuts.tex
\stackMath
\newcommand\tsup[2][2]{%
 \def\useanchorwidth{T}%
  \ifnum#1>1%
    \stackon[-1pt]{\tsup[\numexpr#1-1\relax]{#2}}{\hspace{1pt}\scriptstyle\sim}%
  \else%
    \stackon[.5pt]{#2}{\hspace{1pt}\scriptstyle\sim}%
  \fi%
}

\newcommand{\nc}{\newcommand}

\DeclareMathOperator{\Int}{Int}

\nc{\cO}{\mathcal{O}}

\newcommand{\sgo}{{\mathsf{S}_1(\Ga,\Op)}}

\nc{\ram}{\mathsf{Ramsey}}
\nc{\soo}{\mathsf{S}_1(\cO,\cO)}
\nc{\swl}{\mathsf{S}_1(\Om,\Lambda)}
\nc{\goo}{\gone(\cO,\cO)}
\nc{\gwo}{\gone(\Om,\cO)}
\nc{\gwl}{\gone(\Om,\Lambda)}
\nc{\soox}[1]{\mathsf{S}_1(\cO(#1),\cO(#1))}
\nc{\swlx}[1]{\mathsf{S}_1(\Om(#1),\Lambda(#1))}
\nc{\goox}[1]{\gone(\cO(#1),\cO(#1))}
\nc{\gwox}[1]{\gone(\Om(#1),\cO(#1))}
\nc{\gwlx}[1]{\gone(\Om(#1),\Lambda(#1))}

\nc{\sfinoo}{\mathsf{S}_{\mathrm{fin}}(\cO,\cO)}
\nc{\sfinwl}{\mathsf{S}_{\mathrm{fin}}(\Om,\Lambda)}
\nc{\sfinww}[1]{\mathsf{S}_{\mathrm{fin}}(\Om(#1),\Om(#1))}
\nc{\gfinoo}{\gfin(\cO,\cO)}
\nc{\gfinwo}{\gfin(\Om,\cO)}
\nc{\gfinwl}{\gfin(\Om,\Lambda)}
\nc{\sfinoox}[1]{\mathsf{S}_{\mathrm{fin}}(\cO(#1),\cO(#1))}
\nc{\sfinwlx}[1]{\mathsf{S}_{\mathrm{fin}}(\Om(#1),\Lambda(#1))}
\nc{\sfinwwx}[1]{\mathsf{S}_{\mathrm{fin}}(\Om(#1),\Om(#1))}
\nc{\gfinoox}[1]{\gfin(\cO(#1),\cO(#1))}
\nc{\gfinwox}[1]{\gfin(\Om(#1),\cO(#1))}
\nc{\gfinwlx}[1]{\gfin(\Om(#1),\Lambda(#1))}

\nc{\mc}{\mathcal}
\nc{\thusfar}{\my{--- Edited thus far ---}}
\nc{\lei}{\le^\oo}
\nc{\sqsubs}{\sqsubseteq^*}
\nc{\card}[1]{\left|#1\right|}
\nc{\medcard}[1]{\biggl|\,#1\,\biggr|}
\nc{\smallcard}[1]{|\,#1\,|}
\nc{\bds}{bidirectional $\roth$-scale}
\nc{\bfP}{\mathbf{P}}
\nc{\bfF}{\mathbf{F}}
\nc{\bfQ}{\mathbf{Q}}
\nc{\bbT}{\mathbb{T}}
\nc{\bbZ}{\mathbb{Z}}
\nc{\bbN}{\mathbb{N}}
\nc{\bbC}{\mathbb{C}}
\nc{\beq}{\begin{equation}}\nc{\eeq}{\end{equation}}
\nc{\mbq}{\mb{?}}
\nc{\mb}[1]{{\mbox{\textbf{#1}}}}
\nc{\nop}{$\times$}
\nc{\fbn}{\!\!\fbox{\!\nop\!}\!\!}
\nc{\yup}{\checkmark}
\nc{\forces}{\Vdash}
\nc{\name}[1]{\dot{#1}}
\nc{\tf}{\my{FINISHED THUS FAR}}
\nc{\FU}{Fr\'echet--Urysohn}
\nc{\gs}{$\gamma$~space}
\nc{\Gab}{\Gamma_{\mathrm{B}}}
\nc{\Omb}{\Omega_{\mathrm{B}}}
\nc{\Ga}{\Gamma}
\nc{\Om}{\Omega}
\nc{\smallbinom}[2]{\begin{psmallmatrix} #1\\ #2 \end{psmallmatrix}}
\nc{\bgamma}{\smallbinom{\Om}{\Ga}}

\nc{\productive}[2]{(#1,\allowbreak #2)^\x}
\nc{\prdct}[1]{#1^\x}
\nc{\Sel}{\mathsf{S}}
\nc{\sset}[2]{\{\,#1 : #2\,\}}
\nc{\smb}[1]{{\!\!\mb{#1}\!\!}}
\nc{\medset}[2]{{\biggl\{\,#1 : #2\,\biggr\}}}
\nc{\smallmedset}[2]{{\bigl\{\,#1 : #2\,\bigr\}}}
\nc{\set}[2]{{\left\{\,#1 : #2\,\right\}}}
\nc{\seq}[2]{{\la\, #1 : #2\,\ra}}
\nc{\eseq}[1]{#1_0, \allowbreak #1_1, \allowbreak\dotsc} 

\nc{\eseqint}[3]{#1_{#2}, \allowbreak\dotsc,\allowbreak #1_{#3}} 
\nc{\eseqstart}[2]{#1_{#2},\allowbreak #1_{#2+1},\dotsc } 
\nc{\eprod}[1]{#1_{1}\times \allowbreak#1_{2}\times\dotsb}

\nc{\shortprod}[1]{\prod_{n=1}^\infty{#1}_n}

\nc{\eprodint}[3]{#1_{#2}\times \allowbreak\dotsb\times\allowbreak #1_{#3}}

\nc{\seleseq}[1]{#1_1\in \mathcal{#1}_1, \allowbreak #1_2\in \mathcal{#1}_2, \allowbreak\dotsc}
\nc{\cube}{(\Cantor)^\bbN}
\nc{\Match}{\op{Match}}
\nc{\concat}[1]{\hat{\phantom{a}}\langle #1\rangle}
\nc{\poset}{\mathbb{P}}
\nc{\fn}[1]{{\op{Fn}(#1\times\w,2)}}
\nc{\linadd}{\op{linadd}}
\nc{\nonprod}{\non^\x}
\nc{\alephes}{{\aleph_0}}
\nc{\my}[1]{{\color{red}{#1}}}
\nc{\later}[1]{{\color{green} #1}}
\nc{\BTs}[1]{{\color{green} #1 (BT)}}
\nc{\Cp}{\op{C}_\mathrm{p}}
\nc{\Bp}{\op{B}_p}
\nc{\Pa}[8]{\bibitem{#1} {#2}, \emph{#3}, {#4} \textbf{#5} ({#6}), {#7}--{#8}.}
\nc{\tPa}[5]{\bibitem{#1} {#2}, \emph{#3}, {#4}, to appear.}
\nc{\sPa}[4]{\bibitem{#1} {#2}, \emph{#3}, {#4}, submitted.}
\nc{\Bc}[9]{\bibitem{#1} {#2}, \emph{#3}, in: \textbf{#4} (#5), #6 #7, #8--#9.}
\nc{\fD}{\mathfrak{D}}
\nc{\fX}{\mathfrak{X}}
\nc{\Onbd}{\Op_{\mathrm{nbd}}} 
\nc{\Omnb}{\Om_{\mathrm{nbd}}} 
\nc{\od}{\mathfrak{od}}
\nc{\Setting}[7]{\xymatrix@R=4pt@C=7pt{#1\ar@{-}[r]&#2\ar@{-}[r]&#3\\&#4\ar@{-}[u]\\
#5\ar@{-}[uu]\ar@{-}[r] & #6\ar@{-}[u]\ar@{-}[r] & #7\ar@{-}[uu]}}
\nc{\mx}[1]{\begin{matrix}#1\end{matrix}}
\nc{\plim}{p\txt{-}\lim}
\nc{\Bgp}{{\Z^\bbN}}
\nc{\Cgp}{{{\Z_2}^\bbN}}
\nc{\Cite}[1]{\textbf{[#1]}}
\nc{\Next}[1]{{#1^+}}
\nc{\cFin}{\mathrm{cF}}
\nc{\scsp}{\text{-scale space}}
\nc{\cfn}{\text{cofinal}\ }
\nc{\Con}{\text{Concentrated}}
\nc{\Lind}{\text{Lindel\"of}\,}
\nc{\con}{\text{-Concentrated}}
\nc{\lind}{\text{-Lindel\"of}\,}
\nc{\ctbl}{\text{countably }\allowbreak}
\nc{\Hur}{\text{Hurewicz}}
\nc{\intvl}[2]{{[#1(#2),\allowbreak #1(#2\!+\!1))}}
\nc{\Bdd}{\mathbf{B}}
\nc{\Dfin}{\mathfrak{D}_\mathrm{fin}}
\nc{\grbl}{{\mbox{\textit{\tiny gp}}}}
\nc{\bbP}{\mathbb{P}}
\nc{\bbH}{\mathbb{H}}
\nc{\BOfat}{\B_{\Om_{\mathrm{fat}}}}
\nc{\Bgood}{\B_{\mathrm{good}}}
\nc{\compactN}{\cl{\mathbb{N}}}
\nc{\blocks}[2]{\op{cl}_{#2}(#1)}
\nc{\blocksplus}[2]{\op{cl}^+_{#2}(#1)}
\nc{\arx}[1]{\texttt{http://arxiv.org/math/#1}}
\nc{\bq}{\begin{quote}}
\nc{\eq}{\end{quote}}
\nc{\cl}[1]{\overline{#1}}
\nc{\Cl}[2]{\mathrm{cl}_{#1}(#2)}
\nc{\CH}{the Continuum Hypothesis}
\nc{\MA}{Martin's Axiom}
\nc{\Bfat}{\B_\mathrm{fat}}
\nc{\inv}{^{-1}}
\nc{\Cantor}{{2^\w}}
\nc{\bP}{\mathbf{P}}
\nc{\bof}{\op{\fb}}
\nc{\dof}{\op{\fd}}
\nc{\bofF}{\bof(\cF)}
\nc{\sr}[3]{\underset{\mbox{#3}}{\mbox{#1}}}
\nc{\gp}{\binom{\Om}{\Ga}}
\nc{\gpsmall}{\mbox{$\gp$}}
\nc{\gig}{\gimel}
\nc{\gns}{\sone(\Om,\gig)}
\nc{\nsr}[2]{#1}
\nc{\Srg}{{\mathbb{S}}}
\nc{\Srgs}{{\mathbb{S}^*}}
\nc{\NN}{{\w^{\w}}}
\nc{\ZN}{{\Z^{\bbN}}}
\nc{\NNup}{{\bbN^{\uparrow\bbN}}}
\nc{\NNupb}{{b^{\uparrow\bbN}}}
\nc{\Pof}{\op{P}}
\nc{\PN}{{\Pof(\w)}}
\nc{\rothx}[1]{{[#1]^{\mbox{\tiny $\infty$}}}}
\nc{\tx}{{\tilde{x}}}
\nc{\roth}{{[\w]^{\w}}}
\nc{\roths}{{[b]^{\mbox{\tiny $\infty$}}}} 
\nc{\Fin}{\mathrm{Fin}}
\nc{\ici}{[\bbN]^{ \infty, \infty}}
\nc{\Inc}{{\compactN^{\uparrow\bbN}}}
\nc{\powInc}[1]{{\big(\Inc\big)^{#1}}}
\nc{\powFin}[1]{{\big(\Fin\big)^{#1}}}
\nc{\powPN}[1]{{\big(\PN\big)^{#1}}}
\nc{\NcompactN}{{\compactN^\bbN}}
\nc{\Uarrow}{\smash{\big\uparrow}}
\nc{\LE}{\preccurlyeq}
\nc{\GE}{\succcurlyeq}
\nc{\op}{\operatorname}
\nc{\im}{\op{Im}}
\nc{\Span}{\op{span}}
\nc{\maxfin}{\op{maxfin}}
\nc{\ran}{\op{range}}
\nc{\iso}{\cong}
\nc{\Madd}{{\M}^\star}
\nc{\cI}{\mathcal{I}}
\nc{\cJ}{\mathcal{J}}
\nc{\scrA}{\mathscr{A}}
\nc{\scrB}{\mathscr{B}}
\nc{\scrC}{\mathscr{C}}
\nc{\scrD}{\mathscr{D}}
\nc{\scrF}{\mathscr{F}}
\nc{\scrK}{\mathscr{K}}
\nc{\A}{\D\forall}
\nc{\B}{\mathrm{B}}
\nc{\cB}{\mathcal{B}}
\nc{\cZ}{\mathcal{Z}}
\nc{\bB}{\mathbf{B}}
\nc{\BS}{\mathbf{B}(\mathcal{S})}
\nc{\BF}{\mathbf{B}(\mathcal{F})}
\nc{\BU}{\mathbf{B}(\mathcal{U})}
\nc{\cSp}{\mathcal{S}^+}
\nc{\cFp}{\mathcal{F}^+}
\nc{\cUp}{\mathcal{U}^+}

\nc{\BG}{\B_\Ga}
\nc{\BL}{\B_\Lambda}
\nc{\BT}{\B_\Tau}
\nc{\BTstar}{\B_{\Tau^*}}
\nc{\BO}{\B_\Om}
\nc{\DO}{\cD_\Om}
\nc{\KO}{\cK_\Om}
\nc{\CG}{C_\Ga}
\nc{\CL}{C_\Lambda}
\nc{\CT}{C_\Tau}
\nc{\CTstar}{C_{\Tau^*}}
\nc{\CO}{C_\Om}
\nc{\COgp}{C_{\Om^{\grbl}}}
\nc{\CLgp}{C_{\Lambda^{\grbl}}}
\nc{\BOgp}{\B_{\Om}^{\grbl}}
\nc{\BLgp}{\B_{\Lambda^{\grbl}}}
\nc{\sfC}{\mathsf{C}}
\nc{\sfD}{\mathsf{D}}
\nc{\bD}{\mathbf{D}}
\nc{\Tau}{\mathrm{T}}
\nc{\cA}{\mathcal{A}}
\nc{\cK}{\mathcal{K}}
\nc{\cD}{\mathcal{D}}
\nc{\cF}{\mathcal{F}}
\nc{\cS}{\mathcal{S}}
\nc{\cT}{\mathcal{T}}
\nc{\cG}{\mathcal{G}}
\nc{\cY}{\mathcal{Y}}
\nc{\J}{\mathcal{J}}
\nc{\cL}{\mathcal{L}}
\nc{\cM}{\mathcal{M}}
\nc{\cN}{\mathcal{N}}
\nc{\cH}{\mathcal{H}}

\nc{\Op}{\mathrm{O}}
\nc{\rmA}{\mathrm{A}}
\nc{\rmF}{\mathrm{F}}
\nc{\rmB}{\mathrm{B}}
\nc{\rmD}{\mathrm{D}}
\nc{\rmP}{\mathrm{P}}
\nc{\cC}{\mathcal{C}}
\nc{\cP}{\mathcal{P}}
\nc{\bbQ}{\mathbb{Q}}
\nc{\bbR}{\mathbb{R}}
\nc{\bbS}{\mathbb{S}}
\nc{\cU}{\mathcal{U}}
\nc{\cQ}{\mathcal{Q}}
\nc{\Un}{\bigcup}
\nc{\cV}{\mathcal{V}}
\nc{\cR}{\mathcal{R}}
\nc{\tcR}{\tilde{\mathcal{R}}}
\nc{\cW}{\mathcal{W}}
\nc{\Z}{{\mathbb Z}}
\nc{\Impl}{\Rightarrow}
\long\def\forget#1\forgotten{\marginpar{\textcolor{green}{Forgetting...}}}
\nc{\ft}{\mathfrak{t}}
\nc{\fb}{\mathfrak{b}}
\nc{\fc}{\mathfrak{c}}
\nc{\fd}{\mathfrak{d}}
\nc{\fg}{\mathfrak{g}}
\nc{\oo}{\infty}
\nc{\fr}{\mathfrak{r}}
\nc{\fk}{\mathfrak{k}}
\nc{\bidi}{\mathfrak{bidi}}
\nc{\fu}{\mathfrak{u}}
\nc{\fh}{\mathfrak{h}}
\nc{\fp}{\mathfrak{p}}
\nc{\fj}{\mathfrak{j}}
\nc{\fs}{\mathfrak{s}}
\nc{\w}{\omega}
\nc{\x}{\times}
\nc{\Iff}{\Leftrightarrow}

\nc{\nin}{\notin}
\nc{\cat}{\hat{\ }}
\nc{\sub}{\subseteq}
\nc{\spst}{\supseteq}
\nc{\sm}{\setminus}
\nc{\as}{\subseteq^*}
\nc{\les}{\le^*}
\nc{\leinf}{\le^{\infty}}
\nc{\leS}{\le_S}
\nc{\leF}{\le_F}
\nc{\leU}{\le_U}
\nc{\gew}{\geq^\textrm{w}}
\nc{\rest}{\restriction}

\nc{\la}{\langle}
\nc{\ra}{\rangle}
\nc{\E}{\exists}
\nc{\dom}{\op{dom}}
\nc{\cov}{\op{cov}}
\nc{\add}{\op{add}}
\nc{\addmen}{\add(\Men{})}
\nc{\cof}{\op{cof}}
\nc{\cf}{\op{cf}}
\nc{\non}{\op{non}}
\nc{\unif}{\op{non}}
\nc{\COV}{\op{COV}}
\nc{\ADD}{\op{ADD}}
\nc{\COF}{\op{COF}}
\nc{\NON}{\op{NON}}
\nc{\supp}{\op{supp}}
\nc{\impl}{\to}
\nc{\Lp}{\mathcal{L_\p}}
\nc{\Wlog}{without loss of generality}
\newtheorem{thm}{Theorem}[section]
\nc{\bthm}{\begin{thm}} \nc{\ethm}{\end{thm}}
\newtheorem{need}[thm]{Need}
\nc{\bneed}{\begin{need}\color{dg}} \nc{\eneed}{\end{need}}
\newtheorem{prop}[thm]{Proposition}
\nc{\bprp}{\begin{prop}} \nc{\eprp}{\end{prop}}
\newtheorem{fact}[thm]{Fact}
\nc{\bfct}{\begin{fact}} \nc{\efct}{\end{fact}}
\newtheorem{prob}[thm]{Problem}
\nc{\bprb}{\begin{prob}} \nc{\eprb}{\end{prob}}
\newtheorem{lem}[thm]{Lemma}
\nc{\blem}{\begin{lem}} \nc{\elem}{\end{lem}}
\newtheorem{app}[thm]{Application}
\nc{\bapp}{\begin{app}} \nc{\eapp}{\end{app}}
\newtheorem{claim}[thm]{Claim}
\nc{\bclm}{\begin{claim}} \nc{\eclm}{\end{claim}}
\newtheorem{cor}[thm]{Corollary}
\nc{\bcor}{\begin{cor}} \nc{\ecor}{\end{cor}}
\newtheorem{conj}[thm]{Conjecture}
\nc{\bcnj}{\begin{conj}} \nc{\ecnj}{\end{conj}}
\theoremstyle{definition}
\newtheorem{defn}[thm]{Definition}
\nc{\bdfn}{\begin{defn}} \nc{\edfn}{\end{defn}}
\newtheorem{obs}[thm]{Observation}
\nc{\bobs}{\begin{obs}} \nc{\eobs}{\end{obs}}
\theoremstyle{remark}
\newtheorem{rem}[thm]{Remark}
\nc{\brem}{\begin{rem}} \nc{\erem}{\end{rem}}
\newtheorem{cnv}[thm]{Convention}
\nc{\bcnv}{\begin{cnv}} \nc{\ecnv}{\end{cnv}}
\newtheorem{exam}[thm]{Example}
\nc{\bexm}{\begin{exam}} \nc{\eexm}{\end{exam}}
\nc{\bpf}{\begin{proof}} \nc{\epf}{\end{proof}
}
\nc{\be}{\begin{enumerate}}
\nc{\ee}{\end{enumerate}}
\nc{\bi}{\begin{itemize}}
\nc{\bimy}{\my{\begin{itemize}}
\nc{\eimy}{\end{itemize}}}
\nc{\itm}{\item}
\nc{\ei}{\end{itemize}}
\nc{\Subsection}[1]{\goodbreak\subsection*{#1}}

\nc{\sone}{\mathsf{S}_1}
\nc{\sfin}{\mathsf{S}_\mathrm{fin}}
\nc{\ufin}{\mathsf{U}_\mathrm{fin}}
\nc{\Split}{\mathsf{Split}}

\nc{\gone}{\mathsf{G}_1}    
\nc{\tgfin}{\tilde{\mathsf{G}}_\mathrm{fin}}
\nc{\gfin}{\mathsf{G}_\mathrm{fin}}

\nc{\men}[1]{\sfin(\Op(#1),\Op(#1))}
\nc{\sch}{\ufin(\cO,\Omega)}
\nc{\rothb}{\text{Rothberger}}
\nc{\pmen}{\sfin(\Omega,\Omega)}
\nc{\Rothb}{\sone(\Op,\Op)}

\nc{\prothb}{\sone(\Omega,\Omega)}
\nc{\tU}{{\tilde{U}}}
\nc{\tF}{{\tilde{F}}}
\nc{\tY}{{\tilde{Y}}}
\nc{\tX}{{\tsup[1]{X}}}
\nc{\dtX}{{\tsup[2]{X}}}
\nc{\dt}[1]{{\tsup[2]{#1}}}
\nc{\td}{{\tilde{d}}}
\nc{\tz}{{\tilde{z}}}
\nc{\cfd}{\cf(\fd)}

\nc{\msep}{\sfin(\cD,\cD)}
\nc{\rsep}{\sone(\cD,\cD)}
\nc{\cft}{\sfin(\Omega_{\mathbf{0}},\Omega_{\mathbf{0}})}
\nc{\scft}{\sone(\Omega_{\mathbf{0}},\Omega_{\mathbf{0}})}

\nc{\Umen}{U\text{-Menger}}
\nc{\hur}{\ufin(\cO,\Gamma)}
\nc{\tUmen}{\tU\text{-Menger}}

\nc{\Men}{\text{Menger}}
\nc{\Sch}{\text{Scheepers}}

\nc{\aspst}{\prescript{*}{}{\spst}\ }
\nc{\eqs}{=^*}
\nc{\ctblOm}{\Omega_{\mathrm{ctbl}}}
\nc{\GNga}{{\smallbinom{\Om}{\Ga}}}
\nc{\ctblga}{\smallbinom{\ctblOm}{\Ga}}

\nc{\nadd}{\cN_{\mathrm{add}}}
\nc{\ball}{\mathrm{B}}
\nc{\cOunif}{\cO^{\textrm{unif}}}

\nc{\sep}{
\vspace{2cm}
\noindent
\begin{minipage}{\textwidth}
	\textcolor{red}{\rule{\textwidth}{1pt}}
\end{minipage}
}
\nc{\FS}{\op{FS}}
\nc{\sums}{\op{SS}}
\nc{\SG}{\op{SG}}
\nc{\tSG}{\op{\widetilde{SG}}}
\nc{\G}{\op{G}}
\nc{\pBM}{\op{wgM}}
\nc{\FSG}{\op{FSG}}
\nc{\M}{\op{M}}
\nc{\FP}{\op{FP}}
\nc{\nonNadd}{\non(\nadd)}
\nc{\borga}{\Ga_\mathrm{Bor}}
\nc{\pick}{x}
\nc{\gen}{y}
\nc{\nullzind}{\sone(\{\Op_n^{\mathsf{unif}}\}_{n\in\bbN},\Ga)}
\nc{\nullzindf}[1]{\sone(\{\Op_{#1}^{\mathsf{unif}}\}_{n\in\bbN},\Ga)}

\definecolor{dg}{RGB}{42,101,24}

\nc{\myb}[1]{\textcolor{blue}{#1}}
\nc{\mydg}[1]{{\color{dg}{#1}}}
\DeclareMathOperator{\eexists}{\exists}
\DeclareMathOperator{\fforall}{\forall}
\nc{\Exists}[1]{\bigl(\eexists #1\bigr)}
\nc{\Forall}[1]{\bigl(\fforall #1\bigr)}
\nc{\End}[1]{\bigl(#1\bigr)}
\nc{\dmo}[2]{\DeclareMathOperator{#1}{#2}}
\dmo{\Asc}{Asc}
\nc{\plusmin}{\wedge}

\nc{\cBsub}{{\cB^{\mbox{\tiny $\sub$}}}}

\nc{\Alice}{{\textsc{Alice}{}}}
\nc{\Bob}{{\textsc{Bob}}}
\nc{\BM}{\op{BM}}
\nc{\Palpha}{{\bbS_\alpha}}
\nc{\Pbeta}{{\bbS_\beta}}
\nc{\Pwtwo}{\bbS_{\w_2}}
\nc{\restrict}{\upharpoonright}
\nc{\U}{\mathcal U}
\nc{\zrost}{\bar{\w}^{\uparrow\w}}
\nc{\pr}{\mathrm{pr}}

\nc{\gM}{\op{gM}}
\nc{\GKO}{\gone(\cK,\cO)}